\documentclass[11pt]{amsart}

\renewcommand{\bar}{\overline}
\newcommand{\eps}{\varepsilon}
\newcommand{\pa}{\partial}
\renewcommand{\phi}{\varphi}

\newcounter{hours}\newcounter{minutes}

\newcommand{\ka}{K\"ahler }

\font\strange=msbm10

\newcommand{\R}{{{\mathchoice  {\hbox{$\textstyle{\text{\strange R}}$}}
{\hbox{$\textstyle{\text{\strange R}}$}}
{\hbox{$\scriptstyle  N\kern-0.3em  R$}}  
{\hbox{$\scriptscriptstyle  R\kern-0.2em  R$}}}}}

\newcommand{\Z}{{{\mathchoice  {\hbox{$\textstyle{\text{\strange Z}}$}}
{\hbox{$\textstyle{\text{\strange Z}}$}}
{\hbox{$\scriptstyle  Z\kern-0.3em  Z$}}
{\hbox{$\scriptscriptstyle  Z\kern-0.2em  Z$}}}}}

\newcommand{\N}{{{\mathchoice  {\hbox{$\textstyle{\text{\strange N}}$}}
{\hbox{$\textstyle{\text{\strange N}}$}}
{\hbox{$\scriptstyle  N\kern-0.3em  N$}}
{\hbox{$\scriptscriptstyle  N\kern-0.2em  N$}}}}}

\usepackage{amsmath,amsthm,amscd}

\title
[]{On the lower bound  estimates of 
sections of the canonical bundles over a Riemann surface}
\author{Zhiqin Lu}
\date{\today}

\address[Zhiqin Lu]
{Department of Mathematics\\
Columbia University\\
New York, NY 10027}
\email[Zhiqin Lu]{lu@cpw.math.columbia.edu}
\thanks{Research supported by NSF grant DMS 9971506}

\newtheorem{theorem}{Theorem}[section]

\newtheorem{lemma}{Lemma}[section]

\newtheorem{prop}{Proposition}[section]

\theoremstyle{remark}
\newtheorem{rem}{Remark}[section]

\begin{document}
\maketitle

\numberwithin{equation}{section}

\tableofcontents

\section{Introductions}

Suppose $M$ is an $n$-dimensional K\"ahler manifold and $L$ is an ample
line bundle over $M$. Let the K\"ahler form of $M$ be $\omega_g$ and the
Hermitian metric of $L$ be $H$. We assume that $\omega_g$ is the curvature
of $H$, that is, $\omega_g=Ric(H)$. 
The K\"ahler metric of
$\omega_g$ is called a polarized 
\ka metric on $M$.

Using $H$ and $\omega_g$,
for any positive integer $m$,
 $H^0(M,L^m)$ becomes a Hermitian inner product space. 
We use the following notations: suppose that $S,T\in H^0(M,L^m)$. Let
$<S,T>_{H^m}$ be the pointwise inner product and
\[
(S,T)=\int_M<S,T>_{H^m}\frac{\omega^n_g}{n!}
\]
be the inner product of $H^0(M,L^m)$. Let
\[
||S||=\sqrt{<S,S>_{H^m}
}
\]
be the pointwise norm. In particular, $||S||(x)$ denotes the pointwise
norm at $x\in M$. Let
\[
||S||_{L^2}=\sqrt{(S,S)}
\]
be the $L^2$-norm of $S$. 

Let $\{S_1,\cdots,S_{d(m)}\}$ be an orthonormal basis of $H^0(M,L^m)$.
The quantity(see~\cite{T5})
\begin{equation}\label{quan}
\sum_{i=1}^{d(m)}||S_i||^2
\end{equation}
plays an important rule in K\"ahler-Einstein geometry and stability of
complex manifolds. 
If $m$ is sufficiently large, then there is a natural embedding
$\phi_m: M\rightarrow CP^{d(m)-1}$ by
$x\mapsto [S_1(x),\cdots,S_{d(m)}(x)]$. 
This follows from the Kodaira's embedding theorem.
The metric
$\frac 1m\phi^*\omega_{FS}$ on $M$ is called the Bergman metric. 
The following formula is an easy but key observation made by
Tian~\cite{T5}
in his proof of the covergence of the Bergman metrics:
\[
\frac 1m\pa\bar\pa\log\sum_{i=1}^{d(m)}||S_i||^2=\frac 1m
\phi^*\omega_{FS}-\omega_g.
\]

A lot of work has been done by several authors 
~\cite{T5},
~\cite{ru},~\cite{sz}, ~\cite{cat} and~\cite{Lu10}
on the estimates
of~\eqref{quan}. However, these works are concentrated on a
single manifold. On the other hand, in studying the stability and
K\"ahler-Einstein
geometry
of  manifolds, we need to study the behavior of~\eqref{quan} for
a
family of manifolds. Tian's work~\cite{T6} on the Calabi Conjecture shows
the
importance and non-triviality of the problem of giving lower bound
estimate 
of~\eqref{quan}
for a family of complex surfaces. 
In the  $n$-dimensional case, Tian~\cite{T9} proved that
~\eqref{quan} has
a positive lower bound depends on the
dimension $n$, the upper bound of the Betti numbers, the positve upper and
lower bound of the Ricci curvature and the $L^n$ norm of the Riemannian
sectional curvature.
In general, 
~\eqref{quan}  maybe a useful tool in studying algebraic fiberations over
a compact K\"ahler manifold.
The paper of interest to us are~\cite{MR83h:14025}
and~\cite{MR86j:14032}.

In this paper, we shall study the behavior of~\eqref{quan} on
Riemann surfaces. Even in the case of Riemann surfaces, the problem
of finding a uniform lower bound of~\eqref{quan} is nontrivial.
In fact, by
the counterexample in \S 3, we know that there is no uniform lower
bound in general.  We proved a partial uniform estimate in \S 4
which, I believe, is the ``right'' one in the sense that it gives all the
information on  stability of the Riemann surfaces.

We also consider the coordinate ring of Riemann surfaces. For smooth
Riemann surfaces $M$ of genus $g\geq 2$, it is well  known that its
coordinate ring is
finitely
generated. That is, there is a positive integer $m_0$, such that
for any $S\in H^0(M, K_M^m)$ with $m>m_0$, we can find $U_i\in H^0(M,
K_M^{m_0}), i\in
I$ and $T_i\in H^0(M, K_M^{m-m_0}), i\in 
I$
such that
\[
S=\sum_{i\in I} U_iT_i,
\]
where $I$ is a finite set. In studying the behavior of 
Riemann surfaces near the boundary of the Teichm\"uller space, we need
some uniform estimates. In this paper, we give a uniform estimate
which will give us the information on  singular Riemann surfaces.  

The above setting is similar to that in the corona problem
in complex analysis.
The corona problem on the unit disk was studied by Carleson
in ~\cite
{MR25:5186}. Carleson's result
stimulates many ideas which proved to be useful for other problems.
An extensive discussion of the Carleson's corona theorem can be found
in~\cite{JBG}. 

Modifying Wolff's~\cite{WT} proof of Carleson's
theorem together with the $\bar\pa$-estimate, we give
a uniform corona estimate in the last section of this paper
as an application of Theorem~\ref{third}.
In order to obtain the result, we  take
 special care
to the points where the injective radius are small.

The organization of the paper is as follows: in \S 2, we give a lower
bound of~\eqref{quan} in terms of the genus $g$ and the injective radius
$\delta$ of $M$. In \S 3, we give a counterexample which shows that the
lower bound must depend on $\delta$. In \S 4, we give the partial uniform
estimate. That is, a lower bound of ~\eqref{quan} at $x\in M$ depending
only on the injective radius $\delta_x$ of $x$. In \S 5, we solve the
uniform corona problem by the partial uniform estimate.

The main results of this paper are the following:

\begin{theorem}\label{main}
Let $M$ be a Riemann surface of the genus $g\geq 2$. Let $K_M$ be the
canonical line  bundle of $M$ endowed with a Hermitian metric $H$. Let 
the curvature $\omega_g$ of $H$ be positive. $\omega_g$ gives a \ka
metric of $M$. Let the curvature $K$ of $\omega_g$ satisfy
\[
-K_1\leq K\leq K_2
\]
for nonnegative constants $K_1, K_2\geq 0$ and let $\delta'$ be the 
injective radius of $M$. Let
\[
\delta=\min(\delta',\frac{1}{\sqrt{K_1+K_2}}).
\]
Then there is an absolute constant $C>0$ such that for $m\geq 2$,
\[
\sum_{i=1}^{d(m)}||S_i||^2\geq  e^{-\frac{Cg^3}{\delta^6}},
\]
where $\{S_1,\cdots,S_{d(m)}\}$ is an orthnormal basis of $H^0(M,K_M^m)$.
\end{theorem}

\begin{theorem}\label{second}
For any $\eps>0$ and $m\geq 2$, there is a Riemann surface
$M$
 of genus
$g\geq 2$ with the constant Gauss curvature $(-1)$ such that
\[
\inf_{x\in M}\sum_{i=1}^{d(m)}||S_i||^2\leq \eps.
\]
\end{theorem}

The theorem disproves the conjecture that
the absolute lower bound exists.

When the injective radius of $M$ goes to zero, the first eigenvalue 
and the Sobolev constant will also go to zero. In this case, Theorem
~\ref{main} gives no information. In the following theorem, we proved that
~\eqref{quan} has a lower bound which depends only on the local
information and is independent to the injective radius of $M$. 
For this reason, we call the result partial uniform estimate. 

\begin{theorem}\label{third}
Let $M$ be a Riemann surface of genus $g\geq 2$ and constant curvature
$(-1)$. Then there are absolute constants $m_0>0$ and $D>0$ such that for
any $m>m_0$
and any $x_0\in M$, there is a section $S\in H^0(M,K_M^m)$ with
$||S||_{L^2}=1$ such that
\begin{equation}
||S||(x_0)\geq\frac{\sqrt m}{D(1+\frac{1}{\sqrt m\delta_{x_0}^2}
e^{\frac{\pi}{\delta_{x_0}}})},
\end{equation}
where $\delta_{x_0}$ is the injective radius of $\delta_{x_0}$.
\end{theorem}

On the coordinate ring $\oplus_{j=0}^{\infty}
H^0(M, K_M^j)$, we have the following

\begin{theorem}
Let $M$ be a Riemann surface
as above.
Then there is an $m_0>0$ such that for any $m>m_0$ and $S\in
H^0(M,K_M^m)$,
there is a decomposition
\[
S=\sum_{i=1}^d S_i
\]
of $S_i\in H^0(M,K_M^m) (i=1,\cdots, d)$
such that
\begin{align}
\begin{split}
&||S_i||_{L^2}\leq C(m,m_0,g)||S||_{L^2}\\
&||S_i||_{L^\infty}\leq C(m,m_0,g)||S||_{L^\infty}
\end{split}
\end{align}
for $i=1,\cdots d$,
and
\[
S_i=T_iU_i
\]
for a basis $U_1,\cdots,U_d$ of $H^0(M,K_M^{m_0})$
and $T_1,\cdots,T_d\in H^0(M,K_M^{m-m_0})$.
\end{theorem}

{\bf Acknowledgment.} The author thanks Professor Tian for
suggesting me this problem and
 his help during 
the preparing of this paper. He also thanks Professor Phong for the
support
and some suggestions of the work.

\section{A lower bound estimate}
Suppose that $M$ is a Riemann surface of genus $g\geq 2$.
Let $K_M$ be the canonical line  bundle over $M$ with the Hermitian
metric $H$. We assume that the curvature $\omega_g$ of $H$ is positive
and $\omega_g$ defines a K\"ahler metric of $M$.

 Let $K$ be
the Gauss curvature of the metric $\omega_g$. Let $K_1$ and $K_2$ be two
nonnegative
constants such that
\begin{equation}\label{kbd}
-K_1\leq K\leq K_2.
\end{equation}
Let $\delta'$ be the injective radius of $M$ with  and
\begin{equation}\label{delta}
\delta=\min(\delta',\frac{1}{\sqrt{K_1+K_2}}).
\end{equation}
 Let $x_0\in M$ be a
fixed
point. Let $U$ be the open set 
\[
U=\{dist(x,x_0)<\delta\}. 
\]
It is well known that at each point of $U$ there is an isothermal
coordinate. 
In the first part of this section,
we prove that there is a holomorphic
function $z$ on $U$ which gives the isothermal coordinate of $U$ with 
the \textsl{required} estimate.

Consider the equation
\begin{equation}\label{defh}
\left\{
\begin{array}{l}
\Delta h=K/2\\
h|_{\pa U}=0,
\end{array}
\right.
\end{equation}
where $\Delta$ is the (complex) Laplacian of $M$. The solution $h$ exists
and is
unique.
Let $ds^2$ be the Riemann metric of $U$. Then we have

\begin{lemma}\label{lem21}
The metric $e^h ds^2$ on $U$ is a flat metric.
\end{lemma}

{\bf Proof.} A straightforward computation using~\eqref{defh}.

\qed

Since $U$ is an open set which is differmorphic to an open set in the
Euclidean plane, we can assume that there are global frames on $U$.
Let $\omega^1$ and $\omega^2$ be 1-forms on $U$ such that 
\[
e^hds^2=\omega_1^2+\omega_2^2.
\]
Let $\omega_{12}$ be the connection 1-form defined by
\begin{equation}\label{stru}
\begin{array}{l}
d\omega_1=\omega_{12}\wedge\omega_2,\\
d\omega_2=-\omega_{12}
\wedge\omega_1.
\end{array}
\end{equation}
Then by Lemma~\ref{lem21}, $d\omega_{12}=0$.
It follows that there is a  real smooth function $\sigma$ on $U$ such that
\begin{equation}\label{conn}
\omega_{12}=d\sigma.
\end{equation}
Let 
\[
\xi=e^{i\sigma}(\omega_1+i\omega_2).
\]
Then by~\eqref{stru} and ~\eqref{conn}, we have
\[
d\xi=0.
\]
Thus  there is a function $z$ on $U$ such that
\[
\xi=dz,
\]
and 
\begin{equation}\label{flat}
e^hds^2=dzd\bar z.
\end{equation}
Either $z$ or $\bar
z$ will be holomorphic because it defines a conformal structure of $U$.
Without losing generality, we assume that $z$ is holomorphic
and at $x_0$, $z=0$.

We have the following lemma:

\begin{lemma}\label{esti}
 Let $\rho$ be the distance to the point
$x_0$. $\rho(x)=dist(x,x_0)$. Then
\begin{equation}\label{27}
\frac 13\rho\leq |z|\leq 3\rho 
\end{equation}
for $\rho<\delta$.
\end{lemma}

{\bf Proof.} By the Gauss Lemma~\cite[page 8]{CE}, the Riemann metric
$ds^2$ can be
written as
\[
ds^2=d\rho^2+f^2(\rho,\theta)d\theta^2
\]
for the polar coordinate $(\rho,\theta)$
where $f(\rho,\theta)$ is a smooth function satisfying
\[
f(0,\theta)=0, \frac{\pa f}{\pa\rho}(0,\theta)=1,
\]
and 
\[
\frac{\pa^2 f}{\pa\rho^2}=-Kf.
\]
By the Hessian comparison theorem~\cite[page 4]{SY}, we have
\[
\Delta\rho\geq\frac{\sqrt{K_2}}{4}\cot\sqrt{K_2}\rho.
\]
In particular, $\Delta\rho\geq 0$ on $U$.
Noting that $\Delta$ is the complex Laplacian, we have 
\begin{equation}\label{dist}
\Delta\rho^2=\frac 12|\nabla\rho |^2+2\rho\Delta\rho\geq \frac 12.
\end{equation}

By ~\eqref{dist}, we have
\[
\begin{array}{l}
\Delta(h+K_1\rho^2)\geq K/2+K_1/2\geq 0,\\
\Delta(h-K_2\rho^2)\leq K/2-K_2/2\leq 0.
\end{array}
\]
By the maximal principle, we have
\begin{equation}\label{h}
-1\leq -K_2\delta^2\leq h|_{\pa U}-K_2\rho^2\leq h\leq
h|_{\pa U}+K_1\rho^2\leq K_1\delta^2
\leq 1.
\end{equation}
Let $ds_1^2=e^h ds^2$ denotes the flat metric. Then
\[
e^{-1} ds^2\leq ds_1^2\leq e ds^2.
\]
By~\eqref{flat}, $|z|$ is the distance to the point $x_0$
with respect to the metric $ds_1^2$.
Thus by~\eqref{h},
\[
\frac 13\rho\leq e^{-1}\rho\leq |z|\leq e\rho
\leq 3\rho.
\]

\qed

{\bf Proof of Theorem~\ref{main}.}
Define a smooth function $\eta: \R^+\rightarrow\R$ such that
\begin{equation}\label{eta}
\eta(t)=\left\{
\begin{array}{ll}
0 & t\geq 1,\\
1 & 0\leq t\leq \frac 12.
\end{array}
\right.
\end{equation}
We assume that 
$|\eta'|\leq 4$ and $|\eta''|\leq 4$.

In the rest of this paper $C_1, C_2,\cdots,$ are absolute constants, 
unless otherwise stated.

Let
$\delta_1=\frac 14\delta$.
Define the smooth  function $r$ on $M$ such that 
\begin{equation}\label{deff}
r=\left\{
\begin{array}{ll}
\eta(\frac{|z|}{\delta_1})\log (\frac{|z|}{\delta_1}) &
x\in U,\\
0 & x\notin U.
\end{array}
\right.
\end{equation}
$r$ is well defined.  
For if $x\in\pa U$, then $\rho=\delta$. By the Lemma~\ref{esti},
$|z|\geq\frac 43 \delta_1$ and thus $r|_{\pa U}=0$ using either
expression.

Note that if $z\neq 0$, then $\Delta\log (|z|)=0$. Define the
function
$\psi$ such that $\psi=\Delta r$ for $z\neq 0$ and $\psi=0$ for 
$z=0$. We have

\begin{lemma}\label{lem23}
There is a constant $C_1>0$ such that
\begin{equation}\label{esti2}
\left\{
\begin{array}{l}
|\psi|\leq \frac{C_1}{\delta^2},\\
\int_M|\psi|\leq C_1.
\end{array}
\right.
\end{equation}
\end{lemma}

{\bf Proof.}
A straightforward computation gives
\begin{equation}\label{psi}
\psi=\Delta r=\frac 14 e^h(\frac{1}{\delta_1^2}\eta''\log(
\frac{|z|}{\delta_1})+\eta'\frac{1}{\delta_1|z|}\log(
\frac{|z|}{\delta_1})+2\eta'\frac{1}{\delta_1|z|})
\end{equation}
for $\frac{\delta_1}{2}<z<\delta_1$. Using~\eqref{h},~\eqref{eta},
we have the estimate
\[
|\psi|\leq \frac{C_2}{\delta^2}
\]
for some constant $C_2$.
To get the estimate of the $\int_M|\psi|$, we first see that by the 
volume comparison theorem~\cite[page 11]{SY}, 
\[
vol(U)\leq2\pi(\frac{\cosh \sqrt{K_1}\delta-1}{K_1}).
\]
Since $\sqrt{K_1}\delta\leq 1$, there is a constant $C_3$
such that
\begin{equation}\label{vol}
vol(U)\leq C_3\delta^2.
\end{equation}
The lemma follows by setting $C_1=\max(C_2, C_2C_3)$.

\qed

Let $G(x,y)$ be the Green's function of $M$. That is,
\[
\left\{
\begin{array}{l}
\Delta_x G(x,y)=\frac 14({-\delta_x(y)+\frac{1}{vol(M)}}),\\
\int_M G(x,y) dx=0,
\end{array}
\right.
\]
where $\Delta_x$ is the (complex) Laplacian with respect to $x$
and $\delta_x(\cdot)$ is the Dirac function.
Let $b$ be the function on $M$ such that
\[
\left\{
\begin{array}{l}
\Delta b=K/4+1/2,\\
\int_M b=0.
\end{array}
\right.
\]
Since the K\"ahler metric $\omega_g\in -c_1(M)$, the above 
equation
 has a
unique solution.

Let the function $a: M\rightarrow\R$ defined by
\begin{equation}\label{xxa}
a=G(x,x_0)+\frac{1}{2\pi}(r-\frac{1}{vol(M)}\int_Mr)+b,
\end{equation}
where $x_0$ is the fixed point of $M$ and $r$ is defined in
~\eqref{psi}. Then $a$ is a smooth function on $M$.
We have
\begin{equation}\label{defa}
\left\{
\begin{array}{l}
\Delta a=\frac{1}{4vol(M)}+\frac{1}{2\pi}\psi+K/4+1/2,\\
\int_M a=0.
\end{array}
\right.
\end{equation}

\begin{lemma}\label{lema}
There is a constant $C_5$ such that
\[
|a|\leq\frac{C_5g^3}{\delta^6},
\]
where $g$ is the genus of the Riemann surface $M$.
\end{lemma}

{\bf Proof.} Let $\lambda_1$ be the first eigenvalue of $M$, then 
by the Poincare inequality,
we have
\begin{equation}\label{lambda}
\lambda_1\int_M a^2\leq\int_M|\nabla a|^2.
\end{equation}
Integration by parts using~\eqref{defa}, we have
\begin{equation}\label{nabla}
\int_M|\nabla a|^2\leq\int_M |a(\frac{1}{2\pi}\psi+\frac{1}{4vol(M)}
+\frac{
K_1+K_2}{4}+\frac 12)|.
\end{equation}
Let $g$ be the genus of $M$.
By the Gauss-Bonnet Theorem, $vol(M)=4\pi(g-1)$. On the other hand,
$K_1+K_2\leq\frac{1}{\delta^2}$. Let $a(x')=\max |a|$.
By~\eqref{lambda}, ~\eqref{nabla} and Lemma~\ref{lem23},
\begin{equation}\label{aest}
\int_M a^2\leq\frac{(C_1+6\pi)g}{\lambda_1\delta^2}a(x').
\end{equation}
Consider a neighborhood $U'$ of $x'$ defined by 
\[
U'=\{ x| dist (x,x')	<\delta\}.
\]
Let $z$ be the holomorphic function in Lemma~\ref{esti} such that
$z(x')=0$. Let
\[
U_1=\{|z|<\frac 13\delta\}.
\]
Let $\tilde\Delta=\frac{\pa^2}{\pa z\pa\bar z}$ be the Euclidean Laplacian
on $U_1$. Then by~\eqref{h}, ~\eqref{esti2} and~\eqref{defa}, we have
\begin{equation}
|\tilde\Delta a|\leq 3(2+\frac{C_1+1}{\delta^2}).
\end{equation}
It follows from an elementary fact that there is a constant $C_4$
such that
\begin{equation}\label{aest1}
a(x')\leq C_4(\log\frac 1\delta
+\frac{1}{\delta}
(\int_{U_1} a^2(x))^{\frac 12}).
\end{equation}
On the other hand, Cheeger's inequality~\cite[page 91]{SY} gives
\begin{equation}\label{chee}
\lambda_1\geq\frac{1}{4g^2}\delta^2.
\end{equation}
Combining~\eqref{aest},
~\eqref{aest1} and ~\eqref{chee}, we have the
required estimate.

\qed

Let 
\begin{equation}\label{phi1}
\phi=-4\pi(G(x,x_0)+b).
\end{equation}
Then for $x\neq x_0$, we have
\begin{equation}\label{pine}
\frac{\sqrt{-1}}{2\pi}\pa\bar\pa\phi\geq
(-\frac{1}{2vol(M)}-\frac K2-1)\omega_g.
\end{equation}

\begin{lemma}\label{lem25}
There is a constant $C_6$ such that
\[
\phi\leq\frac{C_6g^{3}}{\delta^6},
\]
and for $|z|\leq\delta_1$,
\[
\phi\geq-\frac{C_6g^{3}}{\delta^{6}}+2\log |z|.
\]
\end{lemma}

{\bf Proof.}
By~\eqref{xxa}, 
\[
\phi=-4\pi(a-\frac{1}{2\pi}(r-\frac{1}{vol(M)}\int_M r))
\]
The lemma follows from Lemma
~\ref{lema} and~\eqref{deff}.

\qed

We need the following proposition from Demailly (see~\cite{T5}):

\begin{prop}\label{del}
Suppose that $(M,g)$ is a complete K\"ahler manifold of complex dimension
$n$, $L$ is a line bundle on $M$ with the Hermitian metric $h$,and
$\phi$ is a function on $M$, which can be approximated by a decreasing
sequence of smooth functions
$\{\phi_l\}_{1\leq l<+\infty}$. If
\[
<\pa\bar\pa\phi_l+\frac{2\pi}{\sqrt{-1}}(Ric(h)+Ric(g)),v\wedge\bar v>_g
\geq C||v||_g^2
\]
for any tangent vector $v$ of type $(1,0)$ at any point of $M$
and for each $l$, where $C>0$ is a constant independent of $l$, 
and $<\cdot ,\cdot >_g$ is the inner product induced by $g$, then 
for any $C^\infty$ $L$-valued $(0,1)$-form $u_1$ on $M$ with
$\bar\pa u_1=0$ and $\int_M||u_1||^2e^{-\phi} dV_g$ finite, there 
exists a $C^\infty$ $L$-valued function $u$ on $M$ such that
$\bar\pa u=u_1$ and
\[
\int_M||u||^2 e^{-\phi} dV_g\leq\frac 1C\int_M||u_1||^2
e^{-\phi} dV_g,
\]
where $dV_g$ is the volume form $g$ and the norm $||\cdot||$ is
induced by $h$. The function $\phi$ is called the weight function.
\end{prop}

\qed

Let $L=K_M^m$ for $m\geq 2$. $H^m$ gives the positive Hermitian
metric on $K_M^m$. Let $\omega_g$ be the K\"ahler
form defined by the curvature of $H$ on $K_M$.
Let $\phi_l=\max(\phi,-l)$ for $l\in\Z^+$, where $\phi$ is defined in
~\eqref
{phi1}. Then by ~\eqref{pine}, we have

\begin{equation}\label{low}
<\pa\bar\pa\phi_l+\frac{2\pi}{\sqrt{-1}}(Ric(H)+Ric(\omega_g)),
v\wedge\bar v>\geq
(m-1-\frac{1}{2vol(M)})||v||^2.
\end{equation}

In order to prove Theorem~\ref{main}, we need to proved that 
for any $m\geq 2$ and $x_0\in M$, there is a section
$S\in H^0(M,K_M^m)$ such that
\[
||S||^2(x_0)/||S||^2_{L^2}\geq e^{-\frac{Cg^{3}}{\delta^6}}.
\]
We will use Proposition~\ref{del} to construct such a section.

Let $e^p$ be the local representation of the metric $H$. That is,
\[
-\frac{\sqrt{-1}}{2\pi}\pa\bar\pa p=\omega_g.
\]

Let 
\[
u_1=\bar\pa\rho(\frac{\sqrt m|z|}{\delta})
e^{-m\frac{\pa p}{\pa\bar z}(x_0)\cdot z} (dz)^m.
\] 
Then $u_1\in \Gamma(M,K_M^m)$. By~\eqref{flat},~\eqref{27}
and~\eqref{eta},
we have
\begin{equation}\label{u-1}
||u_1||^2\leq\frac{12m}{\delta^2}e^{m(p-2 Re\frac{\pa p}{\pa\bar z}
(x_0)z)}
\end{equation}
for $\frac{\delta}{2\sqrt{m}}\leq |z|\leq\frac{\delta}{\sqrt m}$. 
Let 
$U_1=\{x||z|<\frac 13\delta\}$.
Let $\tilde\Delta$ be the Euclidean Laplacian, by
~\eqref{h}, we see
that
\[
|\tilde\Delta p|\leq 3
\]
on $U_1$.
Using the Possion formula we see that there is a constant $C_7$ such that
\[
|\tilde\nabla^2 p(z)|\leq\frac{C_7}{\delta^2}
\]
for $|z|\leq\frac 12\delta$. Thus
\[
|p-2 Re\frac{\pa p}{\pa\bar z}(x_0) z
-p(x_0)|\leq\frac{C_7}{m}
\]
for $|z|<\frac{\delta}{m}$.
Using this estimate and ~\eqref{u-1}, we have
\[
||u_1||^2\leq\frac{12m}{\delta^2} e^{C_7} e^{mp(x_0)}
\]
for $\frac{\delta}{2\sqrt{m}}\leq |z|\leq\frac{\delta}{\sqrt m}$.
Thus by Lemma ~\ref{lem25} and~\eqref{vol},
\[
\int_M||u_1||^2 e^{-\phi}=\int_{\frac{\delta}{2\sqrt m}\leq |z|
\leq\frac{\delta}{\sqrt m}}||u_1||^2 e^{-\phi}
\leq \frac{48C_3 m}{\delta^2}
e^{C_7}
 e^{\frac{C_6g^{3}}{\delta^6}}
e^{m p(x_0)}.
\]

By Proposition~\ref{del} and ~\eqref{low}, there is a
$u\in\Gamma(M,K_M^m)$ such that
$\bar \pa u=u_1$ and,
\[
\int_M||u||^2 e^{-\phi}\leq\frac{1}{m-1-\frac{1}{2vol(M)}}
\int_M||u_1||^2 e^{-\phi}.
\]
Using Lemma~\ref{lem25} again, for $m\geq 2$, there is a $C_8$
such that
\begin{equation}\label{ue}
\int_M||u||^2\leq\
\frac{C_8}{\delta^2} e^{\frac{2C_6 g^{3}}{\delta^6}}
e^{mp(x_0)}.
\end{equation}

On the other hand, we have
\begin{equation}\label{pe}
\int_M||\rho(\frac{m|z|}{\delta}) 
e^{-m\frac{\pa p}{\pa\bar z_0}(x_0)z}
(dz)^m||^2
\leq e^{C_7}\int_{|z|\leq\frac{\delta}{\sqrt m}} e^{mp(x_0)}
\leq e^{mp(x_0)} C_3 e^{C_7}\frac{\delta^2}{m}.
\end{equation}

Let $S=\rho(\frac{|z|}{\delta_1})
e^{-m\frac{\pa p}{\pa\bar z}(x_0)z}
(dz)^m-u$. Then $\bar\pa S=0$.
Since
$\int_Me^{-\phi}=+\infty$, $u(x_0)=0$. In particular, $S\neq 0$. 
Using~\eqref{ue},~\eqref{pe}, we have
\[
||S||^2(x_0)/||S||^2_{L^2}\geq
1/\left (
2C_3e^{C_7}\frac{\delta^2}{m}+
2\frac{C_8}{\delta^2}e^{\frac{2C_6g^{3}}{\delta^6}}
\right)
\]
Thus for $m\geq 2$, there is
a $C$ such that
\[
||S||/||S||_{L^2}\geq e^{-\frac{Cg^3}{\delta^6}}.
\]
This completes the proof of Theorem~\ref{main}.

\section{A counterexample}\label{sec3}

In the last section, we give a lower bound estimate of~\eqref{quan}
in terms of the injective radius of $M$.
In this section, we  give a counterexample that the 
uniform estimate is not true. More precisely, we are going to
disapprove the following:

{\bf Conjecture.} Let $K_M$ be the canonical line bundle of  a Riemann
surface $M$ of genus $g\geq 2$ and constant Gauss curvature $(-1)$. Then
for $m$ sufficiently large, there is a number $C(m,g)>0$,
depending only on $m$ and  $g$,
 such that for
any
orthonormal basis $S_1,\cdots,S_d$ of $H^0(M,K_M^m)$, we have
\[
\inf \sum_{i=1}^d||S_i||^2\geq C(m,g).
\]
In order to give the counterexample,
we use the following Collar Theorem of Keen~\cite[p264]{Keen}:

\begin{theorem}[Keen]\label{keen}
Consider the region $T$ of $U$, the upper half plane, bounded
by the curve $r=1$, $r=e^l$, $\theta=\theta_0$ and $\theta=\pi-\theta_0$.
Let
$\gamma$ be a closed geodesic on $M$ with length $l$. Then there is a
conformal
isometric mapping $\phi: T\rightarrow M$ such that $\phi(iy)=\gamma$.
The image $\phi(T)$ of $T$ is called a collar.
Then we can choose $\theta_0$ small enough such that
the area of the Collar is at least $\frac{8}{\sqrt 5}$.
\end{theorem}

The following theorem
gives the counterexample and implies Theorem~\ref{second}:

\begin{theorem}\label{counterexample} 
For any $\eps>0$ and $m\geq 2$, there is a Riemann surface $M$ of
constant curvature $(-1)$ and genus $g\geq 2$ such that there is a point
$x_0\in M$ satisfying
\[
||S||(x_0)\leq\eps
\]
for any $S\in H^0(M,K_M^m)$ with $||S||_{L^2}=1$.
\end{theorem}

The idea of the proof is that when the length of a closed geodesic 
line tends to zero, the collar will be longer and longer in order
to keep the area of the collar having a lower bound. Topologically,
a collar is a cylinder. By expanding the functions on the corollary using
the Fourier series, we can find the suitable $x_0$ and the estimates. 
We begin by discussing some elementary properties of a collar.

Let $R>0$ be a large real number. Let $(\rho, \theta)
\in (-R,R)\times\R$.
Let the group $\Z$ acting on the space $(-R,R)\times\R$ by
\[
(n,\rho,\theta)\mapsto(\rho,\theta+n\delta)
\]
for $n\in\Z$,
where $\delta>0$ satisfies
\[
\delta\sinh R=\eps_1(=\frac{8}{\sqrt{5}}).
\]
as in Theorem~\ref{keen}.
Define the metric 
\[
ds^2=d\rho^2+(\cosh\rho)^2 d\theta^2
\]
on $(-R,R)\times\R$ which descends to a metric on 
\[
C=(-R,R)\times\R/\Z.
\]
The curvature of the metric is $(-1)$. Note that on $C$, $\rho$ is a
global function but $\theta$ is only locally defined.

We call $C$  a collar
of parameter $\delta$.

Let
\begin{equation}\label{z1}
\left\{
\begin{array}{l}
x=\theta\\
y=2\arctan e^\rho-\frac\pi 2.
\end{array}
\right.
\end{equation}
Define $z=x+iy$. Clearly $z$ is not a global function of $C$. But it 
defines a complex structure of $C$.

Let 
\begin{equation}\label{w1}
w=e^{\frac{2\pi i}{\delta}(\theta+2i(\arctan e^\rho-\frac\pi 4))}
=e^{\frac{2\pi i}{\delta}z}.
\end{equation}
Then $w$ is a global holomorphic function on $C$. Consequently
\begin{equation}\label{dz}
dz=\frac{\delta}{2\pi i}\frac{dw}{w}
\end{equation}
is a global holomorphic 1-form on $C$.

Let $f$ be a holomorphic function on 
a neighborhood of
$\bar C$. Then $f$ is a period
function on $[-R,R]\times\R$, satisfying
\[
f(\rho,\theta\pm\delta)=f(\rho,\theta).
\]

Let the Fourier expansion of $f(-R,\theta)$ and $f(R,\theta)$ be
\begin{align*}
& f(-R,\theta)=\sum_{k=-\infty}^{+\infty} A_ke^{\frac{2\pi
i}{\delta}k\theta};\\
&f(R,\theta)=\sum_{k=-\infty}^{+\infty}B_k
e^{\frac{2\pi i}{\delta}k\theta}.
\end{align*}
Define
\begin{align}\label{g}
\begin{split}
&g_1=\sum_{k=1}^\infty A_ke^{-\frac{4\pi}{\delta}k(\frac{\pi}{4}-\arctan
e^{-R})}w^k,\\
& g_2=B_0,\\
&g_3=\sum_{k=-\infty}^{-1}B_k e^{-\frac{4\pi}{\delta}k(\frac\pi 4-
\arctan e^R)} w^k,
\end{split}
\end{align}
where $w$ is in~\eqref{w1}.
We have the following lemma:

\begin{lemma}\label{lem31}
With the notations as above,
$g_1,g_2,g_3$ are holomorphic functions on $C$. Furthermore
\[
f=g_1+g_2+g_3.
\]
\end{lemma}

{\bf Proof.} 
$g_2$ is a constant. So it is automatically holomorphic.
By equation ~\eqref{w1}, we have
\[
|w|\leq e^{-\frac{4\pi}{\delta}(\arctan e^\rho-\frac\pi 4)}.
\]
Thus we have
\begin{align*}
&|g_1|= |\sum_{k=1}^\infty A_ke^{-\frac{4\pi}{\delta}k(
\frac{\pi}{4}-\arctan e^{-R})} w^k|\\
&\leq\sum_{k=1}^\infty |A_k|e^{-\frac{4\pi}{\delta}k
(\arctan e^\rho-\arctan e^{-R})}\\
&\leq(\sum_{k=1}^\infty|A_k|^2)^{\frac 12}
(\sum_{k=1}^\infty e^{-\frac{8\pi}{\delta}k(\arctan e^\rho-\arctan
e^{-R})})^{\frac 12}.
\end{align*}
By the Bessel inequality,
we have
\[
\sum_{k=1}^{+\infty}|A_k|^2\leq
\frac 1\delta\int_M|f(-R,\theta)|^2 d\theta.
\]
Thus
if $\rho>-R$, the series is convergent absolutely.
So $g_1$ defines a holomorphic function on $\{-R<\rho< R\}$.

By the same argument, $g_3$ is also holomorphic.

\qed

In order to prove that 
\[
f=g_1+g_2+g_3,
\]
we just need to prove that on the set $\{\rho=0\}$, $f=g_1+g_2+g_3$.
Define
\begin{equation}\label{pd}
p_k(\rho)=\int_0^\delta f(\rho,\theta)e^{-\frac{2\pi i}{\delta}
\theta k} d\theta
\end{equation}
for $k\in\Z$.
Apparently
\[
p_k(R)=B_k\delta,\qquad p_k(-R)=A_k\delta
\]
for $k\in\Z$.
By the definition of $z$ in~\eqref{z1}, we have
\[
\frac{\pa}{\pa\bar z}=\frac 12
(\sqrt{-1}\cosh\rho\frac{\pa}{\pa\rho}+\frac{\pa}{\pa\theta}).
\]
Using the equation $\frac{\pa f}{\pa\bar z}=0$, from~\eqref{pd},
we have
\[
p_k'=-\frac{2\pi k}{\delta\cosh\rho}p_k.
\]
Solving the above differential equation gives 
\begin{align}\label{p}
\begin{split}
&p_k(0)=A_ke^{-\frac{4\pi}{\delta}k(\frac\pi 4-\arctan e^{-R})},
\qquad  k=1,2,\cdots,\\
& p_0(0)=B_0,\\
& p_k(0)=B_k e^{-\frac{4\pi}{\delta}k(\frac\pi 4-\arctan e^R)},
\qquad k=-1,-2,\cdots.
\end{split}
\end{align}

From ~\eqref{g} and~\eqref{p}, we see that 
\[
f|_{\rho=0}=(g_1+g_2+g_3)|_{\rho=0}.
\]
Thus
\[
f=g_1+g_2+g_3
\]
on $C$
because both sides are holomorphic functions.

\qed

Let $M$ be a Riemann surface of curvature $(-1)$ and genus $g\geq 2$.
Assume that there is a closed geodesic $\gamma$ on $M$ such that
$length(\gamma)=\delta>0$. Assume that $\delta$ is small enough. Let
$\theta$ be the arc length parameter and $\rho$ be the distance function
to the geodesic. Then

\begin{lemma} $(\rho,\theta)$ is the local coordinate system of $M$ as
long
as 
\[
\sinh\rho\cdot\delta\leq\eps_1=\frac{8}{\sqrt{5}}.
\]
\end{lemma}

{\bf Proof.} 
Note that the area of $\{-R_0<\rho<R_0\}$ is $\delta\sinh R_0$
for any $R_0>0$. The lemma follows from Theorem~\ref{keen}.

\qed

Let $C=\{-R<\rho<R\}$, where
$R$ satisfies $\delta\sinh R=\eps_1$.
Then $z=x+iy$
defines a complex structure of $C$
where $x$ and $y$ are in~\eqref{z1}. We have

\begin{lemma}
Either $z=x+iy$ or $z=x-iy$ is holomorphic on $M$.
\end{lemma}

{\bf Proof.}
A straightforward computation gives
\[
ds^2=(\cosh\rho)^2 dz d\bar z.
\]
Thus $z$ defines a conformal structure which is the same as the one on
$M$. So either $z$ or $\bar z$ is holomorphic.

\qed
 
Without losing generality, we assume that $z$ is
holomorphic.
Fixing $m\geq 2$. 
Let $S\in H^0(M,K_M^m)$. We choose an $x_0\in M$ as follows:
let $\rho_0<0$ be the number such that
\[
\delta (\cosh\rho_0)^2=\eps_1.
\]
Then $\rho_0\rightarrow -\infty$, $\rho_0+R\rightarrow +\infty$ as
$\delta\rightarrow 0$.

By~\eqref{dz}, we see that $(dz)^m\in H^0(C,K_M^m)$. Furthermore
$(dz)^m\neq 0$ on $C$. Thus for any $S\in H^0(M,K_M^m)$, there is
a holomorphic function $f$ on $C$ such that
\[
S|_C=f (dz)^m.
\]
Let
\[
f=g_1+g_2+g_3,
\]
where $g_1,g_2,g_3$ are defined in Lemma~\ref{lem31}. Let
\begin{equation}\label{defs}
S_i=g_i (dz)^m,\qquad i=1,2,3.
\end{equation}

\begin{lemma} With the notations as above, let
$x_0=(\rho_0,0)$. Then
\[
\underset{\delta\rightarrow 0}{\lim}\frac{||S_i||^2(x_0)}
{||S_i||^2_{L^2(C)}}=0
\]
for $i=1,2,3$.
\end{lemma}

{\bf Proof.} By~\eqref{w1}, we have 
\[
w(x_0)=e^{-\frac{4\pi}{\delta}(\arctan e^{\rho_0}-\frac\pi 4)}.
\]
Thus
\begin{align}\label{s1}
\begin{split}
&||S_1||^2(x_0)=\frac{1}{(\cosh\rho_0)^{2m}}
|\sum_{k=1}^\infty A_ke^{-\frac{4\pi}{\delta}k(\arctan e^{\rho_0}
-\arctan e^{-R})+\frac{2\pi i}{\delta}k\theta}|^2\\
&\leq
\frac{1}{(\cosh\rho_0)^{2m}}
\sum_{k=1}^\infty |A_k|^2e^{-\frac{8\pi}{\delta}k(\arctan
e^{\rho_0-1}
-\arctan e^{-R})}\\
\qquad\qquad\qquad
&\cdot
\sum_{k=1}^\infty e^{-\frac{8\pi}{\delta}k(\arctan e^{\rho_0}-\arctan 
e^{\rho_0-1})}\\
&\leq\sum_{k=1}^\infty|A_k|^2e^{-\frac{8\pi}{\delta}k
(\arctan e^{\rho_0-1}-\arctan e^{-R})}\cdot
\sum_{k=1}^{\infty}e^{-\frac{4\pi k}{\delta\cosh(\rho_0-1)}}.
\end{split}
\end{align}
We assume that $\delta$ is so small that
\[
\frac{4\pi}{\delta\cosh(\rho_0-1)}\geq\frac{\mu}{\sqrt{\delta}}>0,
\]
where $\mu$ is an absolute constant. In addition, assume that
$e^{-\frac{\mu}{\sqrt{\delta}}}<\frac 12$. Then

\begin{equation}\label{s1-1}
\sum_{k=1}^\infty e^{-\frac{4\pi k}{\delta\cosh(\rho_0-1)}}
\leq\sum_{k=1}^\infty e^{-\frac{\mu k}{\sqrt{\delta}}}
\leq 2e^{-\frac{\mu}{\sqrt{\delta}}}.
\end{equation}
On the other hand
\[
||S_1||^2_{L^2(C)}=\delta\int_{-R}^R\frac{1}{
(\cosh\rho)^{2m-1}}\sum_{k=1}^\infty |A_k|^2 e^{-\frac{8\pi k}
{\delta} (\arctan e^\rho-\arctan e^{-R})} d\rho.
\]
Assuming $\rho_0-1+R>1$ and $\cosh R<2\sinh R$, we have
\begin{align}\label{s1-2}
\begin{split}
&||S_1||^2_{L^2(C)}\\
&\geq\delta\int_{-R}^{\rho_0-1}
\frac{1}{(\cosh R)^{2m-1}}\sum_{k=1}^\infty
|A_k|^2 e^{-\frac{8\pi k}{\delta}(\arctan e^{\rho_0-1}
-\arctan e^{-R})} d\rho\\
&\geq \frac{\delta}{(\frac{\eps_1}{2\delta})^{2m-1}}
\sum_{k=1}^\infty |A_k|^2 e^{-\frac{8\pi k}{\delta}
(\arctan e^{\rho_0-1}-\arctan e^{-R})}.
\end{split}
\end{align}

By~\eqref{s1}, ~\eqref{s1-1} and~\eqref{s1-2}, we have
\begin{equation}\label{s1-c}
\frac{||S_1||^2(x_0)}{||S_1||^2_{L^2(C)}}
\leq\frac{2^{2m-1}\eps_1^{2m-1}}{\delta^{2m}e^{-\frac{\mu}{\sqrt{\delta}}}}
\rightarrow 0
\end{equation}
for $\delta\rightarrow 0$.

The idea for estimating $S_2$ and $S_3$ are almost the same. 
By~\eqref{defs}, we have
\begin{align}\label{s2}
\begin{split}
& ||S_2||^2(x_0)=\frac{1}{(\cosh\rho_0)^{2m}}|B_0|^2,\\
&||S_2||^2_{L^2(C)}=\delta\int^R_{-R}
\frac{1}{(\cosh\rho)^{2m-1}}|B_0|^2 d\rho.
\end{split}
\end{align}
If $\delta\rightarrow 0$, then $R\rightarrow +\infty$. Thus if
$R$ is large enough, we have
\[
\int_{-R}^R\frac{1}{(\cosh\rho)^{2m-1}} d\rho\geq\mu_1>0,
\]
where $\mu_1$ is a constant only depending on $m$. Thus
\begin{equation}\label{s2-c}
\frac{||S_2||^2(x_0)}{||S_2||^2_{L^2(C)}}=\frac{1}{\mu_1}
\frac{1}{\delta(\cosh\rho_0)^{2m}}
=\frac{1}{\mu_1\eps_1^m}\delta^{m-1}\rightarrow 0
\end{equation}
for $\delta\rightarrow 0$.

For $S_3$, we have
\begin{align}\label{s3}
\begin{split}
&||S_3||^2(x_0)=\frac{1}{(\cosh\rho_0)^{2m}}
|\sum_{k=-\infty}^{-1}B_ke^{-\frac{4\pi}{\delta}k
(\arctan e^{\rho_0}-\arctan e^R)}|^2\\
&\leq\sum_{k=-\infty}^{-1}|B_k|^2 e^{-\frac{8\pi}{\delta}k(\frac\pi
4-\arctan e^R)}
\sum_{k=-\infty}^{-1} e^{-\frac{8\pi}{\delta}k(\arctan e^{\rho_0}-
\frac\pi 4)}.
\end{split}
\end{align}
Since $\rho_0\rightarrow -\infty$, we can assume
\[
\arctan e^{\rho_0}<\frac{\pi}{8}.
\]
Thus
\begin{equation}\label{s3-1}
\sum_{k=-\infty}^{-1} e^{-\frac{8\pi}{\delta}k(\arctan e^{\rho_0}-\frac\pi
4)}
\leq\sum_{k=-\infty}^{-1}e^{\frac{\pi^2}{\delta}k}
\leq 2e^{-\frac{\pi^2}{\delta}}
\end{equation}
for $e^{-\frac{\pi}{\delta^2}}<\frac 12$.
On the other hand,
\begin{align}\label{s3-2}
\begin{split}
&||S_3||^2_{L^2(C)}=\delta
\int_{-R}^R\frac{1}{(\cosh\rho)^{2m-1}}\sum_{k=-\infty}^{-1}
|B_k|^2 e^{-\frac{8\pi k}{\delta}(\arctan e^\rho-\arctan e^R)} d\rho\\
&\geq \delta\int^R_0\frac{1}{(\cosh\rho)^{2m-1}}
\sum_{k=-\infty}^{-1}|B_k|^2 e^{-\frac{8\pi k}{\delta}(\frac\pi 4-\arctan
e^R)} d\rho\\
&\geq\frac{\delta R}{(\cosh R)^{2m-1}}\sum_{k=-\infty}^{-1}
|B_k|^2 e^{-\frac{8\pi k}{\delta}(\frac\pi 4-\arctan e^R)}.
\end{split}
\end{align}

By ~\eqref{s3},~\eqref{s3-1} and~\eqref{s3-2}, we have
\begin{align}\label{s3-c}
&\frac{||S_3||^2(x_0)}{||S_3||^2_{L^2(C)}}
\leq\frac{2e^{-\frac{\pi^2}{\delta}}}{\delta R}
(\cosh R)^{2m-1}
\leq \frac{2^{2m}}{\delta^{2m}}2e^{-\frac{\pi^2}{\delta}}
\rightarrow 0.
\end{align}

Thus for any $\eps>0$ and $m\geq 2$,
from ~\eqref{s1-c}, ~\eqref{s2-c} and~\eqref{s3-c},
 we can find $M$ 
such that there is a closed geodesic with the length sufficiently small
and
an $x_0\in M$ such that
\[
\frac{||S_i||^2(x_0)}{||S_i||^2_{L^2(C)}}\leq\eps
\]
for $i=1,2,3$.

One can check that 
\[
(S_i, S_j)_{L^2(C)}=0.
\]
Thus for any $S\in H^0(M,K_M^m)$ with $||S||_{L^2(M)}$ and
$x_0\in M$,
\[
||S||^2(x_0)\leq\frac{||S||^2(x_0)}{||S||^2_{L^2(C)}}
\leq 3\frac{\sum_{i=1}^3||S_i||^2(x_0)}
{\sum_{i=1}^3||S_i||^2_{L^2(C)}}\leq 3\eps.
\]
Theorem~\ref{counterexample} is proved.

\qed

\section{Partial uniform estimates}
Let
$M$ be a Riemann surface of genus $g$ and constant
 curvature $(-1)$.
In this section, we
prove that there is a (positive) lower bound of ~\eqref{quan} depending
only
on the injective radius of the point. More precisely,
for $m$ large 
enough, for any $x\in M$, there is a section $S\in H^0(M,K_M^m)$
such that $||S||_{L^2}=1$ and $||S||(x)\geq C(\delta_x)$
where $C(\delta_x)$ is a positive constant
depending only on $\delta_x$
 and $\delta_x$ is the
injective radius  at $x$.

Note that in the result the lower bound 
doesn't depend on the injective radius of $M$, which will go to zero as 
$M$ approaches the boundary of the Teichm\"uller space.

We  use all the notations of in \S ~\ref{sec3} about the
collars and the functions on them. The following proposition is a
corollary of the collar theorem:

\begin{prop}\label{prop31}
Let $M$ be a Riemann surface of genus $g\geq 2$ and constant
curvature $(-1)$. Let $\gamma_1,\cdots,\gamma_s$ be the closed geodesics
on $M$ such that
\[
length(\gamma_i)\leq\frac{1}{1000},\qquad
1\leq i\leq s.
\]
Let $C_{\gamma_i} (1\leq i\leq s)$ be the corresponding collars embedded
in
$M$ (Theorem~\ref{keen}). Then for any 
$x\in M\backslash U_{i=1}^s C_{\gamma_i}$,
there is an absolute constant $\eps_2>0$ such that
\[
\delta_x\geq\eps_2.
\]
\end{prop}

{\bf Proof.} 
 Let $\delta=\frac{1}{4000}$.
Let $x\in M\backslash\cup_{i=1}^sC_{\gamma_i}$ and $inj(x)\geq\delta_x>0$.
Then there
is a
point
$y\in M$
such that there are two geodesics $l_1$ and $l_2$
connecting $x$ and $y$ but $l_1$ and $l_2$ are not
homotopic to each other.
If $\delta_x>1$, the theorem has been proved. Otherwise, let $\gamma'$
be the shortest closed curve homotopic to the closed curve $l_1^{-1}l_2$
in
\[
D'=M\backslash\underset{length(\gamma_i)<2\delta}{\cup}
C_{\gamma_i}(R_i-2),
\]
where 
\[
C_{\gamma_i}(R_i-2)=\{
x| dist(x,\gamma_i)< R_i-2\}.
\]
Since
$D'$ is a compact set. If $\gamma'$ doesn't touch any of the
boundary $\pa C_{\gamma_i}(R_i-2)$ for any $i$, then $\gamma'$ must be a
closed
geodesic and by the definition, we have $length(\gamma)\geq
\frac{1}{1000}$ and
thus
$\delta_x\geq\frac{1}{2000}$. Otherwise either $\delta_x>1$ or
$\gamma'\subset
C_{\gamma_i}(R_i-1)\backslash C_{\gamma_i}(R_i-2)$ for some $i$. In the
latter
 case, since $\gamma'$ is not homotopic to zero,
we see  that
\[
length(\gamma')\geq length(\gamma_i) \cosh(R_i-2)\geq\frac {1}{18}\eps_1
\]
(remember $length(\gamma_i)\sinh R_i=\eps_1$).
Thus
\[
\delta_x\geq\frac 12 length(\gamma')\geq \frac{1}{36}\eps_1\geq\eps_2
\]
for $\eps_2=\frac{1}{36}\eps_1$.

\qed

Using the above lemma, we know that outside the collars whose shortest
closed geodesics are small, the injective radius has a lower bound and
the weight function in Proposition~\ref{del} can be constructed in the
ordinary way.
If $x\in C_{\gamma_i}$ for some $i$, we are going to construct the  weight
functions having the compact support within $C_{\gamma_i}$. For this
reason, let's first assume that $C_\delta$ is a collar with 
$\delta<\frac{1}{1000}$ and do some analysis on it.

Let's fix some notations:
there are
absolute constants $\eps_3, \eps_4>0$ such that
\begin{equation}\label{eps2}
\eps_3<\delta\cosh R, \delta\cosh(R\pm 4),\delta\sinh (R\pm 4),
\delta e^{R\pm 4}<\eps_4.
\end{equation}

Let $(\rho,\theta)$ be the local coordinate of the collar
$C=C_\delta$ as in \S ~\ref{sec3}. Then
\begin{equation}\label{3-defw}
w=e^{\frac{2\pi i}{\delta}\theta-\frac{4\pi}{\delta}(\arctan
e^\rho-\frac\pi 4)}
\end{equation}
is the holomorphic function on $C_\delta$.
Let $x_0$ and $p_0$ be the points on $C_\delta$ such that
the local coordinate of $x_0$ and $p_0$ can be represented as:
 $x_0=(\rho_0,0)$ for $R-4>\rho_0\geq 0$ and
$p_0=(R-1,0)$. The function $w$ at $x_0$ and $p_0$
has the values
\begin{align}\label{3-defw0}
\begin{split}
&w_0=e^{-\frac{4\pi}{\delta}(\arctan e^{\rho_0}-\frac\pi 4)},\\
& w_{p_0}=e^{-\frac{4\pi}{\delta}(\arctan e^{R-1}-\frac\pi 4)}
\end{split}
\end{align}
at $x_0$ and $p_0$ respectively.
Let
\begin{equation}\label{alpha}
\alpha=\frac{2\arctan e^{\rho_0}}{\pi},
\end{equation}
and define the functions $\phi_1$, $\phi_2$ and $\phi_3$ 
on $C_\delta$ 
to be

\begin{align}\label{defphi}
\begin{split}
&\phi_1=\log |\frac{w}{w_0}-1|,\\
&\phi_2=\log |\frac{w}{w_{p_0}}-1|,\\
&\phi_3=\phi_1-\alpha\phi_2.
\end{split}
\end{align}

The Riemann metric on $C_\delta$ can be represented as
\begin{equation}\label{rie}
ds^2=d\rho^2+(\cosh\rho)^2 d\theta^2.
\end{equation}
Let the injective radius at $x_0$, $p_0$ and $x$ be
$\delta_{x_0}$, $\delta_{p_0}$ and $\delta_x$. Then we have an
absolute
constant $\eps_5>0$ such that
\begin{equation}\label{ddelta}
\left\{
\begin{array}{l}
\frac 12\delta\cosh\rho_0>\delta_{x_0}>\eps_5\delta\cosh\rho_0,\\
\frac 12\eps_1>\delta_{p_0}>\eps_5,\\
\frac 12\delta\cosh\rho>\delta_x>\eps_5\delta\cosh\rho.
\end{array}
\right.
\end{equation}

We establish some elementary properties of the function
$\phi_3$. 
Let $d=d(x)$ be the distance function to the point $x_0$. Then we have

\begin{lemma}\label{4-lem31}
With the notations as above,
there are constants $C_9, C_{10}>0$ such that
\begin{equation}\label{3-phi3}
\phi_3\leq C_9,
\end{equation}
for $-R+1\leq\rho\leq R-2$
and 
\begin{equation}\label{3-phi31}
\phi_3\geq\log d(x)-\frac{4\pi}{\delta e^\rho}-C_{10}
\end{equation}
for $d(x)\leq\delta_{x_0}$.
\end{lemma}

{\bf Proof.} By ~\eqref{3-defw} and~\eqref{3-defw0}, we have
\begin{align}\label{wfrac}
\begin{split}
&\frac{w}{w_0}=e^{\frac{2\pi i}{\delta}\theta-\frac{4\pi}{\delta}
(\arctan e^\rho-\arctan e^{\rho_0})},\\
&\frac{w}{w_{p_0}}=e^{\frac{2\pi i}{\delta}\theta-\frac{4\pi}{\delta}
(\arctan e^\rho-\arctan e^{R-1})}.
\end{split}
\end{align}

From ~\eqref{wfrac},
we have
\begin{align}\label{basic}
\left\{
\begin{array}{ll}
 \log |\frac{w}{w_0}-1|\leq\log 2&\rho_0\leq\rho\leq R-2,\\
\log|1-\frac{w_0}{w}|\leq\log 2 &-R+1\leq\rho\leq\rho_0,\\
\log |\frac{w}{w_{p_0}}|\geq 0& -R+1\leq\rho\leq R-1,\\
 \log|1-\frac{w_{p_0}}{w}|
\geq\log(1-e^{-\frac{2\pi}{\eps_4}})&
-R+1\leq\rho\leq R-2.
\end{array}
\right.
\end{align}

If $\rho_0\leq\rho\leq R-2$, then we can write
\[
\phi_3=\log|\frac{w}{w_0}-1|-\alpha\log|\frac{w}{w_{p_0}}|
-\alpha\log|1-\frac{w_{p_0}}{w}|.
\]
By ~\eqref{basic}, we have
\begin{equation}\label{phi3-1}
\phi_3\leq\log 2-\log (1-e^{-\frac{2\pi}{\eps_4}}).
\end{equation}

If $-R+1\leq\rho\leq\rho_0$. We can write
\begin{equation}\label{3-phi3r}
\phi_3=\log|\frac{w}{w_0}|-\alpha\log|\frac{w}{w_{p_0}}|
+\log|1-\frac{w_0}{w}|-\alpha\log|1-\frac{w_{p_0}}{w}|.
\end{equation}
Using~\eqref{wfrac},
we have
\begin{align}\label{mainp}
\begin{split}
&\log|\frac{w}{w_0}|-\alpha\log|\frac{w}{w_{p_0}}|\\
&=-\frac 8\delta (\frac\pi 2-\arctan e^{\rho_0})
\arctan e^\rho
+\frac{8}{\delta}\arctan e^{\rho_0}
(\frac\pi 2-\arctan e^{R-1}).
\end{split}
\end{align}
Thus by~\eqref{eps2} and~\eqref{ddelta}
\begin{equation}\label{mains}
-\frac{4\eps_5}{\delta_{x_0}}\leq
\log|\frac{w}{w_0}|-\alpha\log|\frac{w}{w_{p_0}}|
\leq\frac{4\pi}{\eps_3}.
\end{equation}
By~\eqref{basic},~\eqref{3-phi3r} and ~\eqref{mains}, we have
\begin{equation}\label{phi3-2}
\phi_3\leq\frac{4\pi}{\eps_3}+\log 2-\log (1-e^{-\frac{2\pi}{\eps_4}}).
\end{equation}
Thus by ~\eqref{phi3-1} and ~\eqref{phi3-2}, we have $\phi_3\leq C_9$ for
\[
C_9=\frac{4\pi}{\eps_3}+\log 2-\log (1-e^{-\frac{2\pi}{\eps_4}}).
\]

Let's now assume that $d(x)<\delta_{x_0}$. Then by the triangle
inequality
we have 
\begin{equation}
\left\{
\begin{array}{ll}
|\rho-\rho_0|+|\theta|\cosh\rho_0\geq d(x) & 0\leq\theta
\cosh\rho_0<d(x),\\
|\rho-\rho_0|+|\theta-\delta\cosh\rho_0|\geq d(x) & \delta\cosh\rho_0-d(x)
\\&
<\theta\cosh\rho_0<\delta\cosh\rho_0.
\end{array}
\right.
\end{equation}
Without losing generality we assume that $0\leq\theta\leq d(x)$
and $|\rho-\rho_0|+|\theta|\cosh\rho_0\geq d(x)$. If $\theta
\cosh\rho_0\geq\frac 14 d(x)$,
then
\begin{equation}\label{3-sub-1}
|\frac{w}{w_0}-1|\geq
e^{-\frac{4\pi}{\delta}(\arctan e^\rho-\arctan e^{\rho_0})}
\sin\frac{2\pi}{\delta}\theta
\geq\frac\pi 4 e^{-\pi} d(x).
\end{equation}
On the other hand, if $|\rho-\rho_0|\geq\frac 12 d(x)$, then 
\begin{equation}\label{3-sub-2}
|\frac{w}{w_0}-1|\geq e^{-\pi}\frac{2\pi}{\eps_1}d.
\end{equation}
By~\eqref{3-sub-1} and~\eqref{3-sub-2},
there is a
constant $C_{11}>0$ such that
\begin{equation}\label{3-sub-3}
\log|\frac{w}{w_0}-1|\geq\log d -C_{11}.
\end{equation}
We also have
\begin{equation}\label{3-last}
\log|\frac{w}{w_{p_0}}-1|\leq\log 2+\frac{4\pi}{\delta e^\rho}.
\end{equation}
for $d\leq\delta_{x_0}$. By~\eqref{3-sub-3} and ~\eqref{3-last},
from ~\eqref{defphi}
\[
\phi_3\geq\log d-C_{11}-(\log 2+\frac{4\pi}{\delta e^\rho}).
\]
This completes the proof of Lemma~\ref{4-lem31}.

\qed

\begin{lemma}\label{4-lem32}
There is a constant $C_{12}>0$ such that
\[
|\phi_3|\leq C_{12},\qquad |\nabla\phi_3|\leq C_{12}
\]
for $R-3<|\rho|<R-2$.
\end{lemma}

{\bf Proof.} By the above lemma, we see that
\[
\phi_3\leq C_9
\]
for the  constant $C_9$. 
Thus we just need to prove the lower bound of $\phi_3$ and the 
bound for the derivative of $\phi_3$.

If $R-3<\rho<R-2$, 
by~\eqref{wfrac},
we have
\begin{equation}
|\frac{w}{w_0}|\leq e^{-\frac{2\pi}{\eps_1}}.
\end{equation}
Thus
\begin{equation}
\phi_1\geq\log(1-e^{-\frac{2\pi}{\eps_1}}).
\end{equation}
Also we have
\begin{equation}
|\frac{w}{w_{p_0}}|\leq e^{\frac{4\pi}{\eps_3}}
\end{equation}
for $R-3<\rho<R-2$.
Thus
\begin{equation}\label{3-sub-lower1}
\phi_3\geq\log(1-e^{-\frac{2\pi}{\eps_1}})
-\log(1+e^{\frac{4\pi}{\eps_3}}).
\end{equation}
If $R-3<-\rho<R-2$, then by ~\eqref{mainp}, we have
\[
\log|\frac{w}{w_0}|-\alpha\log|\frac{w}{w_{p_0}}|
\geq-\frac{8}{\delta}(\frac\pi 2-\arctan e^{\rho_0})\arctan e^{-(R-3)}
\geq -\frac{2}{\eps_2}.
\]
By~\eqref{3-phi3r}, we have
\begin{equation}\label{3-sub-lower2}
\phi_3\geq -\frac{8}{\eps_3}+\log (1-e^{-\frac{\pi^2}{2}})
-\log 2.
\end{equation}
Combining ~\eqref{3-sub-lower1} and~\eqref{3-sub-lower2}, we get the 
lower bound of $\phi_3$.
Next let's  consider $\nabla\phi_3$. Obviously
\[
|\nabla\phi_3|\leq|\nabla\phi_1|+|\nabla\phi_2|.
\]
Thus we just need to estimate
$|\nabla\phi_1|$ and $|\nabla\phi_2|$. 
By~\eqref{rie}, the Riemann metric under the
coordinate $w$ can be written as 
\[
ds^2=\frac{\delta^2(\cosh\rho)^2}{4\pi^2|w|^2}dw d\bar w.
\]
Thus
\begin{align*}
&|\nabla\phi_1|^2=\frac{4\pi^2|w|^2}{\delta^2(\cosh\rho)^2}
\cdot\frac{1}{|w-w_0|^2}=\frac{4\pi^2}{\delta^2(\cosh\rho)^2}
\cdot\frac{1}{|1-\frac{w_0}{w}|^2},\\
&|\nabla\phi_2|^2=\frac{4\pi^2|w|^2}{\delta^2(\cosh\rho)^2}
\cdot\frac{1}{|w-w_{p_0}|^2}
=\frac{4\pi^2}{\delta^2(\cosh\rho)^2}     
\cdot\frac{1}{|1-\frac{w_{p_0}}{w}|^2}.
\end{align*}
Using the same elementary estimates as above, we get,
\begin{equation}
\left\{
\begin{array}{ll}
|\frac{w_0}{w}|\geq e^{\frac{2\pi}{\eps_1}}&
R-3<\rho<R-2,\\
|\frac{w_0}{w}|\leq e^{-\frac{\pi^2}{2\delta}}&
R-3<-\rho<R-2,\\
|\frac{w_{p_0}}{w}|\leq e^{-\frac{2\pi}{\eps_1}}&
R-3<\rho<R-2,\\
|\frac{w_{p_0}}{w}|\leq e^{-\frac{\pi^2}{2\delta}}
&
R-3<-\rho<R-2.
\end{array}
\right.
\end{equation}
Using these results, we get the bound for the gradient of $\phi_1$ and
$\phi_2$. 
This completes the proof of the lemma.

 \qed

The following proposition summarizes the technical results of this
section.
\begin{prop}\label{tech}
Suppose $M$ is a compact Riemann surface of genus $g\geq 2$ and
constant curvature $(-1)$. Then for any $x_0\in M$, there is a function
$\phi=\phi_{x_0}$ such that $\phi$ is smooth on $M\backslash\{x\}$ and 
\begin{enumerate}
\item In a neighborhood $U_x$ of $x$, $\phi$ can be written as
\[
\phi=2\log d(x)+\psi,
\]
where $\psi$ is a smooth function on $U_x$. Consequently,
\begin{equation}\label{infty}
\int_{U_x} e^{-\phi}=+\infty.
\end{equation}
\item There is a constant $C_{13}$ such that
\begin{equation}\label{t-lower}
\frac{\sqrt{-1}}{2\pi}\pa\bar\pa\phi\geq -C_{13}\omega_g
\end{equation}
on $M\backslash\{x\}$, where $\omega_g$ is the K\"ahler form
of $M$;
\item $\phi$ satisfies
\begin{equation}\label{t-upper}
\phi\leq C_{13}
\end{equation}
on $M$ and
\begin{equation}\label{t-lowerb}
\phi\geq 2\log d(x)-\frac{2\pi}{\delta_{x_0}}-C_{13}
\end{equation}
for $d(x)\leq\delta_{x_0}$.
\end{enumerate}
\end{prop}

{\bf Proof.} Let $\gamma_1,\cdots,\gamma_s$ be the closed geodesics
such that $length(\gamma_i)<\frac{1}{1000}$. Let
$C_{\gamma_i}$ be the corresponding collars. Let
\[
C_{\gamma_i}(R_i-4)=\{
x| dist(x,\gamma_i)\leq R_i-4\}\quad i=1,\cdots s.
\]

For any $x_0\in M$, if $x_0\in C_{\gamma_i}(R_i-4)$ for some $i$,
then let $C_\delta=C_{\gamma_i}$ and
 define
$\phi=\phi_{x_0}:M\rightarrow \R$ as follows

\begin{equation}\label{defphi1}
\phi=\left\{
\begin{array}{ll}
2\eta(\rho-(R-3))\phi_3 & \rho\geq 0\quad and\quad x\in C_\delta,\\
2\eta(-\rho-(R-3))\phi_3 & \rho<0\quad and\quad x\in C_\delta,\\
0 & otherwise,
\end{array}
\right.
\end{equation}
where the function $\eta$ is defined in ~\eqref{eta}.
By Lemma~\ref{4-lem31}, Lemma~\ref{4-lem32} and the fact that
$\phi_3$ is harmonic on $C_\delta\backslash\{x_0\}\backslash\{p_0\}$,
it is easy to check that the function $\phi$ satisfies all the  
assertions in the proposition. On the other hand, if $x_0\notin
C_{\gamma_i}(R_i-4)$ for any
$1\leq i\leq s$, then by Proposition~\ref{prop31}, $\delta_{x_0}\geq
\eps_2$. The argument becomes quiet standard: define
\begin{equation}\label{defphi2}
\phi=2\eta(\frac{d(x)}{\eps_2})\log(\frac{d(x)}{\eps_2})
\end{equation}
Then we can prove that $\phi$ satisfies all the requirements by using the
same
method as in
Lemma~\ref{lem23}.

\qed

\begin{theorem}
Let $M$ be a Riemann surface of genus $g\geq 2$ and constant curvature
$(-1)$. Then there are absolute constants $m_0>0$ and $D>0$ such that for
any $m>m_0$
and any $x_0\in M$, there is a section $S\in H^0(M,K_M^m)$ with
$||S||_{L^2}=1$ such that
\begin{equation}
||S||(x_0)\geq\frac{\sqrt m}{D(1+\frac{1}{\sqrt m\delta_{x_0}^2}
e^{\frac{\pi}{\delta_{x_0}}})}.
\end{equation}
\end{theorem}

{\bf Proof.}
Let $x_0\in M$ and $U_{x_0}=\{x| dist(x,x_0)<\delta_{x_0}\}$. Let $z_1$ be
the holomorphic function  on $U_{x_0}$ such that the hermitian metric
can be represented as
\[
ds^2=\frac{1}{(1-\frac 14|z_1|^2)^2}dz_1 d\bar z_1.
\]
For $m>0$ large enough, let
\[
u_1=\bar\pa(\eta(\frac{2|z_1|}{\delta_{x_0}}))(dz_1)^m.
\]
Then $\bar \pa u_1=0$ and 
\[
||\bar\pa
u_1||^2\leq\frac{16}{\delta^2_{x_0}}
(1-\frac 14|z_1|^2)^{m+1}
\]
for $\frac 14\delta_{x_0}\leq|z_1|\leq\frac 12\delta_{x_0}$.
Thus
there is a $C_{14}>0$ such that
\[
\int_M||\bar\pa u_1||^2 e^{-\phi_{x_0}}\leq
\frac 1m\frac{C_{14}}{\delta^4_{x_0}} e^{\frac{2\pi}{\delta_{x_0}}}.
\]

Let $m_0=C_{13}+2$. 
By Proposition~\ref{del}, for $m>m_0$,  
 we can find a $u\in\Gamma(M,K_M^m)$ such that $\bar\pa u=u_1$ with
\begin{equation}\label{defu1}
\int_M||u||^2e^{-\phi_{x_0}}\leq
\frac{1}{m(m-C_{13}-1)}
\frac{C_{14}}{\delta^4_{x_0}}e^{\frac{2\pi}
{\delta_{x_0}}}.
\end{equation}
Let $S=\eta(\frac{2|z_1|}{\delta_{x_0}}) (dz_1)^m-u$. Then
$\bar\pa S=0$. Thus $S$ is an element of $H^0(M,K_M^m)$. Furthermore,
since $\int_M e^{-\phi_{x_0}}=+\infty$, we must have $u(x_0)=0$. So
\begin{equation}\label{u0}
||S||(x_0)=1.
\end{equation}
On the other hand,
\[
||S||^2_{L^2}\leq 2(\int_M||u||^2+\int_M||\eta (dz_1)^m||^2).
\]
By~\eqref{defu1} and ~\eqref{t-upper}, we have
\begin{equation}\label{u2}
\int_M||u||^2\leq
\frac{1}{m(m-C_{13}-1)}\cdot
\frac{C_{14} e^{C_{13}}}{\delta^4_{x_0}}
e^{\frac{2\pi}{\delta_{x_0}}}.
\end{equation}
We also have
\begin{equation}\label{u3}
\int_M||\eta (dz_1)^m||^2\leq\frac{\pi}{m}.
\end{equation}
The theorem follows from ~\eqref{u0},~\eqref{u2}
and~\eqref{u3}.

\qed

\section{The uniform corona problem}
Let $M$ be a Riemann surface of genus $g\geq 2$. It is well known 
that the coordinate ring $\oplus_{m=0}^\infty H^0(M,K_M^m)$
is finitely generated. That is, there is an $m_0>0$ such that for any
$m>0$ and $S\in H^0(M,K_M^m)$, $S$ can be represented by
\begin{equation}\label{corona}
S=\sum_{i=1}^{d} U_iT_i,
\end{equation}
where $U_i\in H^0(M,K_M^{m_0})$ and $T_i\in H^0(M,K_M^{m-m_0})$
for $i=1,\cdots, d=\dim H^0(M,K_M^{m_0})$.
Finding a suitable set of $\{T_i\}_{i=1,\cdots,d}$ is called the 
corona problem(cf. ~\cite{JBG}).

We need to consider the case where
$M$ approaches to the boundary of the moduli space
in the Teichm\"uller theory. So
in addition to the existence of $U_i$ and $T_i$, we need some uniform 
estimates. 
In this section, we give the uniform estimate for the corona problem on
Riemann surfaces.

\begin{theorem}\label{5main}
Let $M$ be a Riemann surface 
of genus $g$ and constant curvature $(-1)$.
Then there is an $m_0>0$ such that for any $m>m_0$ and $S\in
H^0(M,K_M^m)$,
there is a decomposition
\[
S=\sum_{i=1}^d S_i
\]
of $S_i\in H^0(M,K_M^m) (i=1,\cdots, d)$
such that
\begin{align}\label{result}
\begin{split}
&||S_i||_{L^2}\leq C(m,m_0,g)||S||_{L^2}\\
&||S_i||_{L^\infty}\leq C(m,m_0,g)||S||_{L^\infty}
\end{split}
\end{align}
for $i=1,\cdots d$,
and 
\[
S_i=T_iU_i
\]
for a basis $U_1,\cdots,U_d$ of $H^0(M,K_M^{m_0})$
and $T_1,\cdots,T_d\in H^0(M,K_M^{m-m_0})$.
\end{theorem}

\begin{rem}
An estimate on $T_i (i=1,\cdots, d)$ alone is not expected because of the
counterexample in \S 3, where $\sum||U_i||^2$ can be arbitrarily small. 
\end{rem}

Throughout this section, we will use the notation $D_1, D_2,\cdots$ to
denote the
constants depending only on $m_0$ and the genus $g$.
We also use $A\leq\sim B$ to mean that there is a
positive constant $C=C(m_0,
g)$,
depending only on $m_0$ and $g$ such that $A\leq CB$. Likewise, we use
$A\geq\sim B$ to denote the fact $A\geq CB$
for some constant $C=C(m_0, g)$.

The idea of the proof is that, if the injective radius of $M$ is greater
than
an absolute  constant, then $\sum||U_i||^2$ has a lower bound by an
absolute
 positive constant. In this case, we can solve the corona problem
exactly using the method in~\cite{JBG}. So  we just need to prove the
theorem in the case where $inj(M)$ is arbitrarily small. By the collar
theorem, we know that in this case, there are finite many collars
$C_{\delta_1},\cdots, C_{\delta_s}$ 
(with $\max\delta_i$ small)
embedded in $M$ and they do not
intersect each other. By Proposition~\ref{prop31}, outside the collars,
the injective radius has an absolute lower bound. Special
care must be taken for the sections over these collars.
In order to take care of the collars to get the estimates, 
we first fix a
collar $C_\delta$ embedded in $M$ with the parameter $\delta$
small. We will use all the notations about collars in \S~\ref{sec3}.
For any $\tilde R>0$, 
let
\[
C_\delta(\tilde R)=\{x||\rho|\leq\tilde R\}.
\]
In particular, $C_\delta=C_\delta(R)$ with
$\delta\sinh R=\eps_1$.

We choose and fix a number $m_0>0$ such that $K_M^{m_0}$ is very ample. 
Let $\tilde\eta$ be the cut-off function of $M$ defined as
\begin{equation}\label{etai}
\tilde\eta=
\left\{
\begin{array}{ll}
\eta(\rho-(R-1)) &\rho\geq 0\quad and\quad x\in C_\delta,\\
\eta(-\rho-(R-1)) & \rho<0\quad and\quad x\in C_\delta,\\
0 & otherwise,
\end{array}
\right.
\end{equation}
where the function $\eta$ is defined in ~\eqref{eta}.
Let
\begin{equation}\label{5-u1}
u_1=\frac{1}{\sqrt\delta}\tilde \eta(dz)^{m_0}
\end{equation}
be a section of $K_M^{m_0}$ over $M$ using this cut-off function,
where $dz$ is defined
in~\eqref{dz}.
We can check that
\[
\left\{
\begin{array}{ll}
||\bar\pa u_1||^2\leq\frac{4}{\delta(\cosh\rho)^{2m_0}} &
R-1\leq|\rho|\leq
R,\\
\bar\pa u_1=0 & otherwise.
\end{array}
\right.
\]
 Thus
\[
\int_M||\bar\pa u_1||^2\leq 2\int_{R-1}^R
\frac{4}{\delta(\cosh\rho)^{2m_0}}\cosh\rho\delta d\rho
\leq \frac{8}{(\cosh (R-1))^{2{m_0}-1}}.
\]
By Proposition~\ref{del}, there is a section $u$ of $K_M^{m_0}$ such that
\[
\bar\pa u=\bar\pa u_1
\]
with 
\[
\int_M||u||^2\leq\frac{1}{m_0-1}\int_M||\bar\pa u_1||^2
\leq\frac{8}{(m_0-1)(\cosh (R-1))^{2m_0-1}}.
\]
By using ~\eqref{eps2}, we see that
\begin{equation}\label{5-u3}
\int_M ||u||^2\leq\sim\delta^{2m_0-1}.
\end{equation}
Let 
\begin{equation}\label{5-defu}
U'=u_1-u.
\end{equation}
Then $\bar\pa U'=0$.

\begin{lemma}\label{claim}
Let $U''$ be a holomorphic section of $K_M^{m_0}$ 
on $C_\delta(R)$
such that
$(U'',u_1)_{C_\delta(R-2)}=0$. That is,
\begin{equation}\label{prep}
\int_{C_\delta(R-2)} <U'',u_1>=0,
\end{equation}
where $u_1$ is the section defined in~\eqref{5-u1}. Then there
is an absolute constant $\eps_6>0$ such that
\begin{equation}\label{5-n1}
||U''||\leq\sim e^{-\eps_6e^{R-|\rho|}+m_0(R-|\rho|)}||U''||_{L^2
(C_\delta)},
\end{equation}
and
\begin{equation}\label{5-n2}
||\nabla U''||\leq
\sim e^{-\eps_6e^{R-|\rho|}+(m_0+1)(R-|\rho|)}||U''||_{L^2
(C_\delta)}.
\end{equation}
\end{lemma}

We use the notation in \S ~\ref{sec3}. Let
\begin{align*}
& w_{1}=e^{\frac{2\pi\sqrt{-1}}{\delta}\theta
-\frac{4\pi}{\delta}
(\arctan e^{\rho}-\arctan e^{-(R-2)})},\\
& w_{2}=e^{-\frac{2\pi\sqrt{-1}}{\delta}\theta
+\frac{4\pi}{\delta}
(\arctan e^{\rho}-\arctan e^{R-2})}.
\end{align*}
Let   
\[
U''=(g_1(w_1)+a+g_3(w_2))(dz)^{m_0}
\]
be the decomposition similar to that in~\eqref{g} where
$g_1(w_1)$ and $g_2(w_2)$ are holomorphic functions of $w_1$ and
$w_2$ respectively, $g_1(0)=g_2(0)=0$, $a$ is a constant,
and $dz$ is defined in~\eqref{dz}.

Using~\eqref{prep}, we see that $a=0$. By the Schwartz Lemma, we have
\begin{equation}\label{sch}
|g_1(w_1)|\leq |w_1| \max_{|w_1|=1}|g_1(w_1)|.
\end{equation}
At each point of $\{|w_1|=1\}$ or $\{\rho=R-2\}$, by the collar theorem,
there is an
absolute lower bound for the injective radius. Thus by the 
Cauchy integral formula,
we have
\begin{equation}\label{5-8}
\max_{|w_1|=1}||g_1(dz)^{m_0}||\leq
\sim||g_1(dz)^{m_0}||_{L^2(C_\delta)}.
\end{equation}

Using ~\eqref{sch} and~\eqref{5-8}, we have
\begin{equation}\label{5-8-1}
||g_1 (dz)^{m_0}||\leq e^{-\frac{4\pi}{\delta}(\arctan e^\rho-\arctan
e^{-(R-2)})}
\frac{(\cosh R)^{m_0}}{(\cosh\rho)^{m_0}}
||g_1(dz)^{m_0}||_{L^2(C_\delta)}.
\end{equation}
It is elementary to check that there is an absolute constant
$\eps_6>0$ such that
\begin{equation}\label{5-8-1-1}
e^{-\frac{4\pi}{\delta}(\arctan e^\rho-\arctan e^{-(R-2)})}
\leq e^{-\eps_6e^{R-|\rho|}}.
\end{equation}
Combining~\eqref{5-8-1} and~\eqref{5-8-1-1}, 
\begin{equation}\label{5-13a}
||g_1(dz)^{m_0}||\leq
\sim e^{-\eps_6e^{R-|\rho|}+m_0(R-|\rho|)}
||g_1(dz)^{m_0}||_{L^2(C_\delta)}.
\end{equation}
Similarly, we have
\begin{equation}\label{5-14a}
||g_2(dz)^{m_0}||\leq
\sim e^{-\eps_6e^{R-|\rho|}+m_0(R-|\rho|)}
||g_2(dz)^{m_0}||_{L^2(C_\delta)}.
\end{equation}
~\eqref{5-13a} and
~\eqref{5-14a} give the inequality~\eqref{5-n1}.

On the other hand, a straightforward
computations gives
\begin{align}
\begin{split}
&\nabla U''=(\frac{2\pi\sqrt{-1}}{\delta}(w_1 g_1'+w_2
g_2')\\
&\qquad -\frac 12 m_0\sqrt{-1}\sinh\rho (g_1+ g_2))
(dz)^{m_0}\otimes dz.
\end{split}
\end{align}
Using the same argument as above, we get~\eqref{5-n2}.

\qed

\begin{lemma}\label{lem51}
With the notations as above, 
there is a  constant $r>2$, depending only on $m_0$ and the genus
$g$, such
that
 $U'\neq 0$ on $C_\delta(R-r)$, where $U'$ is defined in~\eqref{5-defu}.
\end{lemma}

{\bf Proof.}
By~\eqref{5-defu}, we know that $u$ is holomorphic on $C_\delta(R-2)$. Let
\begin{equation}\label{alpha1}
u=u'+\alpha u_1
\end{equation}
be the decomposition of $u$ such that $(u',u_1)
_{C_\delta(R-2)}=0$ and $\alpha$ is
a constant.
Then
\begin{equation}\label{5-5}
\int_M||u||^2\geq\int_{C_{\delta}(R-2)}
||u'||^2
+\int_{C_{\delta}(R-2)}
|\alpha|^2||u_1||^2.
\end{equation}
In particular
\begin{equation}\label{5-6}
 \int_{C_{\delta}(R-2)}
||u||^2
\geq |\alpha|^2\int^{R-2}_{-(R-2)}\frac{1}{(\cosh\rho)^{2m_0-1}}
d\rho.
\end{equation}
By~\eqref{5-u3} and ~\eqref{5-6}, we have
\begin{equation}\label{5-7}
|\alpha|
\leq\sim \delta^{m_0-\frac 12}.
\end{equation}
On the other hand, by Lemma~\ref{claim},
\begin{equation}\label{5-7-3}
||u'||\leq\sim
e^{-\eps_6e^{R-|\rho|}+m_0(R-|\rho|)}\delta^{m_0-\frac 12}.
\end{equation}
Thus by~\eqref{5-u1}, ~\eqref{5-7} and~\eqref{5-7-3}, there are constants
$D_1$ and $D_2$ such that
\begin{align}\label{5-7-4}
\begin{split}
&||U'(x)||\geq\frac{1}{\sqrt\delta(\cosh\rho)^{m_0}}
(1-D_1\delta^{m_0-\frac 12})\\
&-D_2
e^{-\eps_6e^{R-|\rho|}+m_0(R-|\rho|)}\delta^{m_0-\frac 12}
\end{split}
\end{align}
for $|\rho|<R-3$.
If $r$ is large enough, then $||U'(x)||>0$ for $|\rho|<R-r$.
In particular $U'\neq 0$ on $C_\delta(R-r)$. This completes the proof of
the lemma.

\qed

Now we assume the general case. Let $C_{\delta_1},\cdots, C_{\delta_s}$ be
the collars of parameters $\delta_1,\cdots, \delta_s$ respectively
embedded
into $M$. We assume that $C_{\delta_i}, (i=1,\cdots, s)$ do not intersect
each other.
By Lemma~\ref{lem51}, there are $R_i, (i=1,\cdots, s)$ such that we can
find
$U_i'\in H^0(M,K_m^{m_0}), (i=1,\cdots, s)$ with
$U_i'|_{C_{\delta_i}(R_i)}\neq 0$. Let $\rho_i, dz_i, (i=1,\cdots, s)$ be
defined in \S~\ref{sec3}
coresponding to $C_{\delta_i} (i=1,\cdots, s)$, respectively. Let
\[
\eta_i=
\left\{
\begin{array}{ll}
\eta(\rho_i-(R_i-1)) &\rho_i\geq 0\quad and\quad x\in C_{\delta_i}\\
\eta(-\rho_i-(R_i-1)) &\rho_i\leq 0\quad and\quad x\in C_{\delta_i}\\
0 & otherwise
\end{array}
\right.
\]
be the cut-off functions for $i=1,\cdots, s$. Assume that $\max \delta_i$
is small enough.
Then
by Lemma~\ref{lem51}, $U_1',\cdots, U_s'\in H^0(M,K_M^{m_0})$
have  the following properties 
\begin{enumerate}
\item $||U_i'||\geq\sim\frac{1}{\sqrt{\delta_i}}e^{-m_0|\rho_i|}$ on
$C_{\delta_i}(R_i)$ for $i=1,\cdots, s$ 
(by ~\eqref{5-7-4});
\item
There are decompositions 
 $U_i'|_{C_{\delta_i}(R_i)}=\alpha_iv_{1i}+v_{2i}$ with
\begin{equation}\label{vj}
v_{1i}=\frac{1}{\sqrt{\delta_i}}\eta_i(dz_i)^{m_0},\qquad
1\leq i\leq s
\end{equation}
and
\[
(v_{1i},v_{2i})_{C_{\delta_i}(R_i)}=0, \qquad 1\leq i\leq s,
\]
where $\alpha_i (i=1,\cdots,s)$ are a constant such that for
$1\leq i\leq s$,
\begin{equation}\label{vj1}
1\leq\sim\alpha_i\leq\sim 1;
\end{equation}
\item By~\eqref{5-u1} and~\eqref{5-u3},
\begin{equation}\label{5-7-5}
\int_{M\backslash C_{\delta_i(R_i)}}||U_i'||^2\leq\sim \delta_i^{2m_0-1}
\end{equation}
for $i=1,\cdots, s$.
\end{enumerate}

We have the following lemma:

\begin{lemma}\label{lem531}
With the notations as above,
there are holomorphic sections $U_1,\cdots, U_s\in 
H^0(M,K_M^{m_0})$ such that
\begin{enumerate}
\item $||U_i||\geq\sim \frac{1}{\sqrt\delta_i}e^{-m_0|\rho_i|}$ on 
$C_{\delta_i}(R_i)$;
\item $(U_i|_{C_{\delta_j}(R_j)},v_{1j})_{C_{\delta_j}(R_j)}=0$ for $i\neq
j, 1\leq i,j\leq
s$;
\item $\int_{M\backslash
C_{\delta_i}(R_i)}||U_i||^2\leq\sim
\delta_i^{2m_0-1}$.
\end{enumerate}
\end{lemma}

{\bf Proof.} Let
\[
\beta_{ij}=(U_i',v_{1j})_{C_{\delta_j}(R_j)}, \qquad 1\leq i,j\leq s.
\]
Then if $i\neq j$ we have
\begin{equation}\label{b-1}
|\beta_{ij}|\leq\sim\delta_i^{m_0-\frac 12}
\end{equation}
by~\eqref{5-7-5} and the definition of $v_{1j},
(j=1,\cdots, s)$ in~\eqref{vj}. We also have
\begin{equation}\label{b-2}
1\geq\sim 
\beta_{ii}\geq\sim 1
\end{equation}
by~\eqref{vj1}. 

Let $B=(\beta_{ij})_{s\times s}$ be the matrix of $(\beta_{ij})$ for
$1\leq i,j\leq s$. Then   by~\eqref{b-1}, ~\eqref{b-2}, $B$  is an
invertible matrix,
when $\max\delta_i$ is small enough.
Let $A=B^{-1}$ be the
inverse matrix and let $A=(\alpha_{ij})_{s\times s}$. Define
\[
U_i=\sum_{j=1}^s\alpha_{ij}U_j', \quad 1\leq i\leq s.
\]
Then $U_i (i=1,\cdots, s)$ satisfies all the requirements in the lemma
by the fact that
\[
\begin{array}{ll}
|\alpha_{ij}|\leq\sim \delta_i^{m_0-\frac 12} & i\neq j,\\
1\geq\sim \alpha_{ij}\geq\sim 1 & i=j.
\end{array}
\]

\qed

Let $U_{s+1}',\cdots, U_d'$ be an orthonormal basis of the space 
$span\{ U_1,\cdots. U_s\}^\perp$. Assume that
\[
(U_i', U_j)=0
\]
for $s<i\leq d, 1\leq j\leq s$.
Let
\begin{equation}\label{kk}
U_i=U_i'-\sum_{j=1}^s\gamma_{ij}U_j,\quad s<i\leq d,
\end{equation}
where
\begin{equation}\label{b-3}
\gamma_{ik}=\frac{1}{(U_k, v_{1k})_{C_{\delta_k}(R_k)}}
(U_i', v_{1k})_{C_{\delta_k}(R_k)},
\quad 1\leq k\leq s, s<i\leq d.
\end{equation}
Then we have
\[
(U_i|_{C_{\delta_j}(R_j)}, v_{1j})_{C_{\delta_j}(R_j)}=0
\]
for $s<i\leq d$ and $1\leq j\leq s$.

We have the following lemma:

\begin{lemma}\label{lem53}
If $x\notin \cup_{j=1}^sC_{\delta_j}(R_j-2)$, then 
\[
\sum_{i=1}^d||U_i||^2\geq\sim 1.
\]
\end{lemma}

{\bf Proof.} Let
\[
l_{ij}=(U_i, U_j)
\]
for $1\leq i,j\leq d$. If $1\leq i,j\leq s$, we have
\begin{equation}\label{k1}
\begin{array}{ll}
(U_i, U_j)\leq\sim \max \delta_i & i\neq j,\\
1\geq\sim (U_i, U_j)\geq\sim  1 & i=j
\end{array}
\end{equation}
by~\eqref{b-1}, ~\eqref{b-2} and lemma~\ref{lem531}. 

If $1\leq i\leq s, s<j\leq d$, we have
\[
(U_i, U_j)=-\sum_{k=1}^s\gamma_{jk} (U_i, U_k)
\]
using ~\eqref{kk}.
By the definition of $\gamma_{jk}$ in~\eqref{b-3}, we have
\begin{equation}\label{g1}
|\gamma_{jk}|\leq\sim\max\delta_j.
\end{equation}
Thus
\begin{equation}\label{k2}
|(U_i, U_j)|\leq\sim\max\delta_j
\end{equation}
for $1\leq i\leq s$ and $s<j\leq d$.
Finally, if $s<i,j\leq d$, then 
\begin{equation}\label{k3}
\begin{array}{ll}
|(U_i, U_j)|\leq\sim \max\delta_i & i\neq j,\\
1\geq\sim (U_i, U_j)\geq\sim 1 & i=j
\end{array}
\end{equation}
by~\eqref{kk} and~\eqref{g1}.
Using~\eqref{k1}, ~\eqref{k2} and~\eqref{k3}, we have
\[
\begin{array}{ll}
|l_{ij}|\leq\sim \max \delta_i & i\neq j,\\
1\geq\sim |l_{ij}|\geq\sim 1 & i=j.
\end{array}
\]

Let $(m_{ij})_{d\times d}$ be the matrix such that
\[
\sum_{i=1}^d\sum_{t=1}^d m_{ji}\bar{m_{tk}}l_{it}=\delta_{jk}
\]
for $1\leq j,k\leq d$. We can choose $m_{ij}$ such that
\begin{equation}\label{m-1}
\begin{array}{ll}
|m_{ij}|\leq\sim \max
\delta_i & i\neq j,\\
1\geq\sim m_{ij}\geq\sim 1 & i=j.
\end{array}
\end{equation}

It is easy to check that
$\sum_{j=1}^d m_{ij} U_j$ for $1\leq i\leq d$ forms an orthonormal basis
of $H^0(M, K_M^{m_0})$.
Thus we have
\[
\sum_{i=1}^d||\sum_{j=1}^dm_{ij}U_j||^2\geq\sim 1
\]
by Theorem~\ref{third}.
The lemma thus follows from~\eqref{m-1} and the fact that $\max \delta_i$
is
small.

\qed

We summarize the results up to now in the following
\begin{prop}\label{prop51}
Let $U_1,\cdots, U_d$ be sections of $H^0(M,K_M^{m_0})$
as above. Then
\begin{align}
&||U_i||+||\nabla U_i||\leq\sim 1 &x\notin
\cup_{j=1}^sC_{\delta_j}(R_j-2),
\label{5-88a}\\
&||U_i||+||\nabla U_i||\leq\sim
e^{-(m_0+1)(R_j-|\rho_j|)}
&x\in C_{\delta_j}(R_j-2), 1\leq j\leq s, \label{5-88b}
\end{align}
for $i>s$ and
\begin{align}
&||U_i||+||\nabla U_i||\leq\sim
{\delta_i}^{m_0-\frac 12}
&x\notin\cup_{j=1}^s C_{\delta_j}(R_j-2),\label{5-88c}\\
&||U_i||+||\nabla U_i||\leq\sim
\frac{1}{\sqrt\delta_i e^{m_0|\rho_i|}}
&x\in C_{\delta_i}(R_i-2),\label{5-88d}\\
&||U_i||+||\nabla U_i||\leq\sim\delta_i^{m_0-\frac 12}
e^{-(m_0+1)(R_j-|\rho_j|)} &
x\in C_{\delta_j}(R_j-2), j\neq i.\label{5-88e}
\end{align}
for $1\leq i\leq s$.
Furthermore, 
\begin{equation}\label{new}
\sum_{k=1}^d||U_i||^2\geq\sim 1
\end{equation}
for $x\notin \cup_{i=1}^s C_{\delta_i}(R_i-2)$.
\end{prop}

{\bf Proof.}  If $x\notin \cup_{j=1}^sC_{\delta_j}(R_j)$, then $\delta_x$
has a uniform lower bound.
~\eqref{5-88a} follows from the Cauchy integral formula. ~\eqref{5-88b}
is a corollary of Lemma~\ref{claim}. ~\eqref{5-88c} follows from
~\eqref{5-u3} and the Cauchy integral formula. ~\eqref{5-88d} follows
from
a straightforward computation.~\eqref{5-88e} follows from
Lemma~\ref{claim} and~\eqref{5-7-5}. Finally, ~\eqref{new} is just a
restatement of the conclusion of Lemma~\ref{lem53}.

\qed

Define an inner product $<\,,\,>$ in the coordinate ring $\oplus_{m=0}^
\infty H^0(M,K_M^m)$. Let $S_1\in H^0(M,K_M^m)$ and $S_2\in
H^0(M,K_M^{m_1})$. Suppose that $m\ge m_1$. We define a section of
$K_M^{m-m_1}$ as follows: Let $x\in M$ and $U_x$ is a local trivialization
of 
$K_M$. let $S_1|_{U_x}=S'\cdot S''$ for 
$S'\in \Gamma(U_x, K_M^{m-m_1})$ and $S''\in
\Gamma(U_x, K_M^{m_1})$. Then
\[
S_3|_{U_x}=S'<S'', S_2>_{H^{m_1}},
\]
where $<,\,,\,>_{H^{m_1}}$ is the pointwise inner product.

{\bf Proof of Theorem~\ref{5main}.} 
We modify the method of Wolff's~\cite{WT} of solving the corona problem
on the unit disk.
 First we construct a $C^\infty$ solution. Let $S\in
H^0(M,K_M^m)$ for a fixed $m>m_0$.
Let
\begin{equation}\label{5-b1}
\left\{
\begin{array}{ll}
 b_k=\eta_k\frac{S}{U_k}+(1-\sum_{j=1}^s\eta_j)
\frac{<S, U_k>}{\sum_{j=1}^d||U_j||^2}& 1\leq k\leq s,\\
 b_k=(1-\sum_{j=1}^s\eta_j)\frac{<S, U_k>}{\sum_{j=1}^d||U_j||^2}
& k>s.
\end{array}
\right.
\end{equation}
Here $b_k (1\leq k\leq s)$ is well defined
because of Lemma~\ref{lem51}. We can check that
\[
S=\sum_{k=1}^d U_k b_k.
\]

If $x\notin \cup_{j=1}^sC_{\delta_j}(R_j-2)$, then by Lemma~\ref{lem53},
\begin{equation}\label{5-b2}
\left\{
\begin{array}{ll}
||b_k||\leq\sim
(\sqrt\delta_k (\cosh\rho_k)^{m_0}+1)||S||& 1\leq k\leq s,\\
 ||b_k||\leq\sim ||S||
&  k>s.
\end{array}
\right.
\end{equation}
By~\eqref{5-b1}, we have
\begin{align}\label{5-b3}
\begin{split}
&\bar\pa b_k=\bar\pa\eta_k\frac{S}{U_k}
-\bar\pa\sum_{j=1}^s\eta_j
\frac{<S,U_k>}{\sum_{j=1}^d||U_j||^2}
+(1-\sum_{j=1}^s\eta_j)\frac{<S,\nabla U_k>}
{\sum_{j=1}^d||U_j||^2}\\\quad-
&(1-\sum_{j=1}^s\eta_j)\frac{<S,U_k>\sum_{j=1}^d<U_j,\nabla U_j>}
{(\sum_{j=1}^d||U_j||^2)^2}\qquad 1\leq k\leq s,\\
&\bar\pa b_k=
-\bar\pa\sum_{j=1}^s\eta_j
\frac{<S,U_k>}{\sum_{j=1}^d||U_j||^2}
+(1-\sum_{j=1}^s\eta_j)\frac{<S,\nabla U_k>}
{\sum_{j=1}^d||U_j||^2}\\
&\quad-
(1-\sum_{j=1}^s\eta_j)\frac{<S,U_k>\sum_{j=1}^d<U_j,\nabla U_j>}
{(\sum_{j=1}^d||U_j||^2)^2}\quad  k>s.
\end{split}
\end{align}
Thus
\begin{equation}\label{5-b4}
\left\{
\begin{array}{ll}
||\bar\pa b_k||\leq \sim
\frac{1}{\delta_k^{m_0-\frac 12}}
||S||& 1\leq k\leq s,\\
||\bar\pa b_k||\leq\sim||S|| &  k>s.
\end{array}
\right.
\end{equation}
Let
\begin{equation}
c_{ik}=\frac{<\bar\pa b_k, U_i>}{\sum_{j=1}^d||U_j||^2}.
\end{equation}
for $i\neq k$ and
$1\leq i,k\leq d$. Then 
by~\eqref{5-b2} and~\eqref{5-b4}, we have
\begin{equation}\label{5-c1}
\left\{
\begin{array}{ll}
 ||c_{ik}||\leq\sim \frac{1}{\delta_k^{m_0-\frac 12}}
||S|| & 1\leq k\leq s,\\
 ||c_{ik}||\leq \sim 
||S|| &  k>s.
\end{array}
\right.
\end{equation}

Let's consider the equations
\[
\bar\pa b_{ik}=c_{ik}
\]
for $i\neq k$ and $1\leq i,k\leq d$. By Proposition~\ref{del}, the
solutions exist and we may assume that
\[
||b_{ik}||_{L^2}\leq ||c_{ik}||_{L^2}
\]
for $i\neq k, 1\leq i,k\leq d$. Thus
by~\eqref{5-c1}, we have
\begin{equation}\label{5-b5}
\left\{
\begin{array}{ll}
||b_{ik}||_{L^2}\leq\sim \frac{1}{\delta_k^{m_0-\frac 12}}
||S||_{L^2} & 1\leq k\leq s,\\
||b_{ik}||_{L^2}\leq \sim ||S||_{L^2}& 
 s<k\leq d.
\end{array}
\right.
\end{equation}

Let
\[
T_i=b_i+\sum_{k=1}^d(b_{ik}-b_{ki}) U_k.
\]
One can check that $\bar\pa T_i=0$ and
\[
S=\sum_{i=1}^dT_iU_i.
\]

In order to prove the theorem, we need to estimate $||T_iU_i||$
for $1\leq i\leq d$. By~\eqref{5-b5} and Proposition~\ref{prop51}, we have
\begin{equation}\label{5-e-1}
||T_iU_i||_{L^2(M\backslash \cup_{j=1}^sC_\delta(R_j))}
\leq\sim ||S||_{L^2}.
\end{equation}
 By the Cauchy integral formula and Proposition~\ref{prop31}
\begin{equation}\label{5-e-2}
||T_iU_i||(x)\leq\sim ||S||_{L^2}
\end{equation}
for $1\leq i\leq d$ and $x\notin \cup_{j=1}^sC_{\delta_j}(R_j)$.

Let
\[
E_i(\rho)=\{x\in C_{\delta_i}(R_i)||\rho_i(x)|\geq\rho\}, \quad
1\leq i\leq s
\]
for positive number $\rho>0$.
By~\eqref{5-b5} and Proposition~\ref{prop51} again, we have
\begin{equation}\label{5-b6}
\left\{
\begin{array}{ll}
 ||T_iU_i-\eta_iS||_{L^2(E_j(\rho)-E_j(\rho+2))}\leq
\sim e^{-(R_i-\rho)}||S||_{L^2}& 1\leq i\leq s,\\
 ||T_iU_i||_{L^2(E_j(\rho)-E_j(\rho+2))}\leq\sim e^{-(R_i-\rho)}
||S||_{L^2} &i>s.
\end{array}
\right.
\end{equation}
 By the Cauchy formula,
\begin{equation}\label{5-b7}
\left\{
\begin{array}{l}
||T_iU_i-\eta_iS||(x)\leq\sim \frac{1}{\delta_x}e^{-(R_i-\rho)}
||S||_{L^2}\quad 1\leq i\leq s,\\
||T_iU_i||(x)\leq\sim \frac{1}{\delta_x}e^{-(R_i-\rho)}
||S||_{L^2}\quad i>s.
\end{array}
\right.
\end{equation}
for any $x\in E_j(\rho+\frac 12)-E_j(\rho+1)$.
Since
\[
\delta_x\geq\eps_5\delta_j(\cosh\rho_j)
\]
by~\eqref{ddelta}, we have
\begin{equation}\label{5-b8}
\left\{
\begin{array}{ll}
||T_iU_i-\eta_iS||(x)\leq\sim ||S||_{L^2} & 1\leq i\leq s,\\
||T_iU_i||(x)\leq\sim ||S||_{L^2} & i>s
\end{array}
\right.
\end{equation}
for any $x\in E_j(\rho+\frac 12)-E_j(\rho+1)$ and any $|\rho_j|<R_j-3$,
$j=1,\cdots s$. Thus
if $||S||_{L^2}=1$, then by~\eqref{5-e-2} and~\eqref{5-b8}
\[
||T_iU_i||_{L^2}\leq\sim 1
\]
for $1\leq i\leq d$. If $||S||_{L^\infty}=1$, then
\[
||T_iU_i||_{L^\infty}\leq\sim 1
\]
for $1\leq i\leq d$.
These results give the inequality  ~\eqref{result}.

\qed

\bibliographystyle{abbrv}
\bibliography{bib}

\end{document}